\documentclass[12pt,reqno]{amsart}
\usepackage{amssymb}
\usepackage{bbm}
\usepackage{mathrsfs}
\usepackage{graphicx}
\usepackage{amsfonts}
\usepackage{enumerate}
\usepackage{tikz}
\usepackage{amsmath, latexsym}
\usepackage{hyperref}
\hypersetup{
  hypertexnames=false,
  colorlinks=true,
  linkcolor=blue,
  citecolor=blue,
  urlcolor=blue
}
\usepackage[top=2cm, bottom=2cm, left=1.8cm, right=1.8cm]{geometry}

\newtheorem{thm}{Theorem}[section]
\newtheorem{cor}[thm]{Corollary}
\newtheorem{lem}[thm]{Lemma}
\newtheorem{prop}[thm]{Proposition}
\newtheorem{defn}[thm]{Definition}

\newtheorem{remark}[thm]{Remark}
\numberwithin{equation}{section}

\newcommand{\ee}{\mathbb{E}}
\newcommand{\nn}{\mathbb{N}}
\newcommand{\rr}{\mathbb{R}}
\newcommand{\zz}{\mathbb{Z}}
\newcommand{\T}{\mathbb{T}}

\def\CC{\mathcal C}
\def\EE{\mathcal E}
\def\LL{\mathcal L}
\def\SS{\mathcal S}

\def\bu{\boldsymbol u}
\def\y{\boldsymbol y}
\def\vep{\varepsilon}
\def\<{\langle}
\def\>{\rangle}

\def\bbeq{\begin{equation}}
\def\nneq{\end{equation}}
\def\bdef{\begin{defn}}
\def\ndef{\end{defn}}
\def\bthm{\begin{thm}}
\def\nthm{\end{thm}}
\def\bprop{\begin{prop}}
\def\nprop{\end{prop}}
\def\blem{\begin{lem}}
\def\nlem{\end{lem}}
\def\bcor{\begin{cor}}
\def\ncor{\end{cor}}
\def\bprf{\begin{proof}}
\def\nprf{\end{proof}}

\newcommand{\cij}{\smash[t]{c_{ij}^{\scalebox{0.8}{$\scriptscriptstyle (2)$}}}}
\def\osc{{\bf osc}}
\def\var{{\rm Var}}
\def\cov{{\rm Cov}}
\def\Ent{{\rm Ent}}

\makeatletter
\def\@settitle{\begin{center}%
  \baselineskip14\p@\relax
  \bfseries
  \@title
  \end{center}}
\makeatother

\begin{document}

\date{\today}

\title[Uniform $W_2I$ and log-Sobolev inequalities for Coulomb gases]
{From one particle to $N$ particles: uniform transport-information and log-Sobolev inequalities for Coulomb gases on the torus}

\author[A. Guillin]{A. Guillin}
\address{{\bf Arnaud Guillin}. Universit\'e Clermont Auvergne, CNRS, LMBP, F-63000 Clermont-Ferrand, France.}
\email{arnaud.guillin@uca.fr}

\author[W. Liu]{W. Liu}
\address{{\bf Wei Liu}. School of Mathematics and Statistics, Wuhan University, Wuhan, Hubei 430072, PR China.}
\email{wliu.math@whu.edu.cn}

\author[B. Nectoux]{B. Nectoux}
\address{{\bf Boris Nectoux}. Universit\'e Clermont Auvergne, CNRS, LMBP, F-63000 Clermont-Ferrand, France.}
\email{boris.nectoux@uca.fr}

\author[L. Wu]{L. Wu}
\address{{\bf Liming Wu}. Universit\'e Clermont Auvergne, CNRS, LMBP, F-63000 Clermont-Ferrand, France.}
\email{Li-Ming.Wu@uca.fr}
\maketitle

\vspace{-0.5cm}
\begin{abstract}
We study how a functional inequality satisfied by the one-particle conditional distributions of a Gibbs measure, with a constant independent of the positions of the other particles, can be propagated to the whole $N$-particle system, with a constant independent of $N$, when the interaction is singular. Our model is the Coulomb gas on the torus $\T^3$ in the mean-field regime. We distinguish two mechanisms. The first one is a dependent tensorization, governed by Dobrushin-Wu interdependence coefficients. It yields uniform Poincar\'e and transport-information $W_2I$ inequalities with fully explicit constants. We explain why this mechanism does not extend to the log-Sobolev inequality: the classical criteria of Zegarlinski, Otto-Reznikoff and Marton fail for the Coulomb interaction, however small the charge is, because of its singularity. The second mechanism is a tensorization of the relative entropy with a defect. Using a uniform bound on the partition function valid at every temperature, due to Wang and Zhao, it yields a uniform log-Sobolev inequality for the Coulomb gas on $\T^3$, under the same condition on the charge as the uniform Poincar\'e inequality. We identify the assumptions on the interaction under which each mechanism applies. All the results are then extended  to the two-dimensional Coulomb gas and to periodic Riesz gases $|x|^{-s}$, $0<s\le1$.  The extension to $d \ge 4$ is also discussed.

\end{abstract}

%

\section{Introduction}

\subsection{Notation and functional inequalities}\label{sec:prelim}

\subsubsection{Basic notation}
Throughout this work, we work on the torus $\T^n=\rr^n/\zz^n$, with $n=3$ or $n=3N$ (identifying $(\T^3)^N$ with $\T^{3N}$), endowed with its flat Riemannian metric $d$\footnote{This is the quotient metric inherited from the Euclidean metric of $\rr^n$ via the canonical projection $\pi:\rr^n\to\T^n$. The geodesic distance between $x,y\in\T^n$ is $d(x,y)=\min_{m\in\zz^n}|x-y+m|$. Identifying $\T^n$ with $(-\frac12,\frac12]^n$, it can be computed componentwise as $d(x,y)=\big(\sum_{l=1}^n\min(|x_l-y_l|,\,1-|x_l-y_l|)^2\big)^{1/2}$.}, and with the Lebesgue (Haar) probability measure $dx$. Writing $x=(x^1,\dots,x^N)$ with $x^i\in\T^3$, the metric of $(\T^3)^N$ is $d_N(x,y)= (\sum_{i=1}^N d(x^i,y^i)^2 )^{1/2}$. In Section~\ref{sec:ext2} we also work on $\T^2$, with the same conventions.

Let $\mu$ be a probability measure on $\T^n$. The associated Dirichlet form is given by $\EE_\mu[f]:=\int_{\T^n}|\nabla f|^2\,d\mu$ for smooth $f$ (its closure and its domain $\mathbb D(\EE_\mu)$ are discussed in Section~\ref{subsec:FI} for the measures considered in this paper). Let $\mathcal P_1(\T^n)$ denote the space of probability measures on $\T^n$. For $k\geqslant1$, the Wasserstein-$k$ distance between $\nu,\mu\in\mathcal P_1(\T^n)$ is
$$
W_k(\nu,\mu) := \left(\inf_{P\in\CC(\nu,\mu)}\iint_{\T^n\times\T^n} d(x,y)^k\,P(dx,dy)\right)^{1/k},
$$
where $\CC(\nu,\mu)$ is the set of all couplings of $\nu$ and $\mu$. It is finite since $\T^n$ is compact.

We write $|v|$ for the Euclidean norm of $v\in\rr^n$ and $v\cdot w$ for the Euclidean inner product. Functions on $\T^n$ are identified with $\zz^n$-periodic functions on $\rr^n$. On $(\T^3)^N$, $\nabla_i$ and $\Delta_i$ denote the gradient and the Laplacian with respect to $x^i$, and $|\nabla f|^2=\sum_{i=1}^N|\nabla_if|^2$. For a measure $\mu$, we write $\mu(f)=\int f\,d\mu$, $\<f,g\>_\mu=\int fg\,d\mu$, and $L^p(\mu)$ for the corresponding Lebesgue spaces. We denote the norm of $L^p(\T^n,dx)$ by $\|\cdot\|_p$. Finally, $\nu\ll\mu$ means that $\nu$ is absolutely continuous with respect to $\mu$, and $\osc(f):=\sup_{x,y}(f(x)-f(y))$.

\subsubsection{Entropy and Fisher-Donsker-Varadhan information}
The relative entropy of $\nu\in\mathcal P_1(\T^n)$ with respect to $\mu$ is
$$
H(\nu|\mu) :=
\begin{cases}
\displaystyle\int_{\T^n} \log\!\left(\frac{d\nu}{d\mu}\right)d\nu, & \text{if } \nu\ll\mu,\\[6pt]
+\infty, & \text{otherwise,}
\end{cases}
$$
and, for $g\geqslant0$, $\Ent_\mu(g):=\int g\log g\,d\mu-\mu(g)\log\mu(g)$, so that $H(\nu|\mu)=\Ent_\mu(\frac{d\nu}{d\mu})$. The Fisher-Donsker-Varadhan information of $\nu$ with respect to $\mu$ is
$$
I(\nu|\mu) :=
\begin{cases}
\EE_\mu[h]=\displaystyle\int_{\T^n} |\nabla h|^2\,d\mu, & \text{if } \nu=h^2\mu \text{ with } h=\sqrt{\tfrac{d\nu}{d\mu}}\in\mathbb D(\EE_\mu),\\[6pt]
+\infty, & \text{otherwise.}
\end{cases}
$$
For a smooth positive density $f=\frac{d\nu}{d\mu}$, $I(\nu|\mu)=\frac14\int\frac{|\nabla f|^2}{f}\,d\mu$, which is one fourth of the usual Fisher information.

\subsubsection{Functional inequalities and hierarchy}
We recall four classical functional inequalities for a probability measure $\mu$ on $\T^n$, ordered from strongest to weakest.
The {\it log-Sobolev inequality} (LSI) holds with constant $C_{\mathrm{LS}}>0$ if for all smooth $f$,
$$
\Ent_\mu(f^2) \le  2C_{\mathrm{LS}}\,\EE_\mu[f],
$$
or equivalently $H(\nu|\mu)\le  2C_{\mathrm{LS}}\, I(\nu|\mu)$ for all $\nu\in\mathcal P_1(\T^n)$.
The {\it transport-information $W_2I$ inequality} holds with constant $C>0$ if for all $\nu\in\mathcal P_1(\T^n)$,
$$
W_2^2(\nu,\mu) \le  (2C)^2\,I(\nu|\mu).
$$
The {\it transport-entropy $W_2H$ inequality} (Talagrand) holds with constant $C_T>0$ if for all $\nu\in\mathcal P_1(\T^n)$,
$$
W_2^2(\nu,\mu) \le  2C_T\,H(\nu|\mu).
$$
The {\it Poincar\'e inequality} (PI) holds with constant $C_P>0$ if for all smooth $f$,
$$
\var_\mu(f) \le  C_P\,\EE_\mu[f].
$$
The classical hierarchy of implications is
$$
\text{LSI} \implies W_2I \implies W_2H \implies \text{PI}.
$$
The implication $\text{LSI}\Rightarrow W_2H$ with $C_T=C_{\mathrm{LS}}$ is due to Otto and Villani~\cite{OV00}. Combined with the LSI in the form $H\le  2C_{\mathrm{LS}} I$, it gives $W_2^2\le  2C_{\mathrm{LS}}H\le  (2C_{\mathrm{LS}})^2 I$, that is, $W_2I$ with $C=C_{\mathrm{LS}}$. The implication $W_2I\Rightarrow W_2H$ is established in~\cite{GLWY_PTRF09, GLWW09}, where transport-information inequalities were introduced and compared with the other functional inequalities. Finally, $W_2H\Rightarrow\text{PI}$ with $C_P\le  C_T$ follows by linearization; see~\cite{OV00, BGL} and also~\cite{CGW_PTRF10}.

\subsection{The Coulomb gas on the torus}\label{subsec:torus}
The periodic Coulomb Green function $W$ on $\T^3$ is defined by
\begin{equation}\label{eq:W}
 -\Delta W=4\pi(\delta_0-1),\qquad \int_{\T^3} W\, dx =0.
\end{equation}
Equivalently, its Fourier coefficients are $\hat W(k)= 1/({\pi|k|^2})$ for $k\in\zz^3\setminus\{0\}$ and $\hat W(0)=0$. It is well known that $W(x)=|x|^{-1}+R(x)$ near the origin, with $R$ smooth there, see e.g. \cite{Serfaty26}. In particular, $W$ is bounded from below, and we set $W_*:=\min_{\T^3}W<0$. The Coulomb gas of $N\geqslant 2$ particles on $\T^3$ is the Boltzmann distribution
\bbeq\label{Coulomb-measure}
\mu_N(dx^1,\cdots,dx^N)=\frac{1}{Z(N)}\exp\Big(-\frac{\Gamma}{2(N-1)}\sum_{{i,j=1,  i\neq j}}^{N}W(x^i-x^j)\Big)\,dx^1\cdots dx^N,
\nneq
where $\Gamma>0$ is the charge parameter and $Z(N)$ is the normalization constant. Since the particles move in the compact torus, no confining potential is needed. The factor $\frac1{N-1}$ is the mean-field scaling, also called the diffusive temperature scale: the total charge felt by one particle is of order one. The measure $\mu_N$ is the reversible invariant measure of the Langevin dynamics
\bbeq\label{Langevin_D}
dX_t^i = - \frac{\Gamma}{N-1} \sum_{j=1,j\neq i}^N\nabla W(X_t^i-X_t^j)\, dt + \sqrt{2}\, dB_t^i, \qquad i=1,\cdots, N,
\nneq
valued in $(\T^3)^N$, where $B^1,\cdots,B^N$ are independent standard Brownian motions in $\rr^3$ (since the drift is $\zz^3$-periodic, \eqref{Langevin_D} is understood modulo $\zz^3$). Its generator
\bbeq\label{Langevin_G}
\LL f= \sum_{i=1}^{N} \left(\Delta_i f -  \frac{\Gamma}{N-1} \sum_{j=1,j\neq i}^N \nabla W(x^i-x^j)\cdot\nabla_i f \right)
\nneq
is symmetric on $L^2(\mu_N)$. The singular interaction causes no difficulty: the diffusion \eqref{Langevin_D} is well defined and a.s.\ never reaches the collision set $\SS_N$ defined below, and its generator is identified with the Dirichlet form of $\mu_N$, see Section~\ref{subsec:FI}.

Consider the collision set
$$\SS_N:=\{x\in (\T^3)^N,\ x^i=x^j \text{ for some }i\neq j \}$$
and the finite energy space
$$
H^1(\mu_N) = \left\{ f\in L^2(\mu_N);\ \nabla f\in L^1_{loc}\big((\T^3)^N\backslash \SS_N,dx\big),\ \EE_{\mu_N}[f]:=\int_{(\T^3)^N\backslash \SS_N} |\nabla f|^2\, d\mu_N<+\infty \right\},
$$
where $\nabla f$ is taken in the distributional sense on the open set $(\T^3)^N\backslash \SS_N$. Since $\SS_N$ is Lebesgue-null, $\EE_{\mu_N}[f]= \int_{(\T^3)^N} |\nabla f|^2\, d\mu_N$. For $\nu\in \mathcal P_1 ((\T^3)^N)$, the Fisher-Donsker-Varadhan information of $\nu$ with respect to $\mu_N$ is
\bbeq\label{FI}
I(\nu|\mu_N)=\begin{cases}
             \EE_{\mu_N}[h] , & \mbox{if } \nu= h^2 \mu_N, \ h = \sqrt{\frac{d\nu}{d\mu_N}} \in H^1(\mu_N),\\
             +\infty, & \mbox{otherwise}.
           \end{cases}
\nneq
Probabilistically, $\nu \mapsto I(\nu|\mu_N)$ is the Donsker-Varadhan rate function governing the large deviations of the empirical measures of the Langevin diffusion \eqref{Langevin_D}.

For $i\in\{1,\dots,N\}$ and $x^{-i}:=(x^j)_{j\neq i}$, the conditional distribution of the $i$-th particle under $\mu_N$ is
\bbeq\label{cond_law}
\mu_i(dx\,|\,x^{-i}) = \frac{1}{Z(x^{-i})}
\exp\left(-\frac{\Gamma}{N-1}  \sum_{j=1,j\neq i}^NW(x-x^j)\right)dx .
\nneq
It is the law of one particle in the field of $N-1$ poles located at the other particles, each of charge $\frac{\Gamma}{N-1}$, so that the total charge of the poles is $\Gamma$, independently of $N$.

\subsection{Objectives: from one particle to $N$ particles}\label{subsec:objectives}
The aim of this paper is to explain how a functional inequality established for the one-particle conditional distributions \eqref{cond_law}, with a constant independent of $N$ and of the positions $x^{-i}$ of the poles, can be propagated to the $N$-particle measure $\mu_N$, with a constant independent of $N$, despite the singularity of the interaction. A central theme is the distinction between two mechanisms.

\emph{Dependent tensorization.} For a product measure, the Poincar\'e, $W_2I$ and log-Sobolev inequalities satisfy the independent tensorization: the constant of the product is the worst constant of the factors. For a weakly dependent measure, one may hope that the dependence only costs a multiplicative factor, which tends to $1$ as the dependence vanishes. The classical way to quantify the dependence is through interdependence coefficients, which measure how much the conditional law $\mu_i(\cdot|x^{-i})$ moves when one coordinate $x^j$ moves. Dobrushin's coefficients in total variation \cite{Dobrushin68,Dobrushin70} and their Wasserstein versions \cite{Wu_AOP06} lead to tensorization statements of the form
$$
C_{\rm global}\le\frac{C_{\rm one\ site}}{1-r},
$$
where $r<1$ is the norm of the interdependence matrix. We show that this mechanism works for the Poincar\'e and $W_2I$ inequalities of the Coulomb gas: the relevant coefficients only involve integrated quantities of the interaction, namely $\|W\|_{L^1}$ and a weighted $H^{-1}$-norm of $\nabla W$, which are finite in dimension $3$.

\emph{Tensorization of relative entropy with a defect.} For the log-Sobolev inequality, the criteria used in the literature do not apply to the Coulomb kernel, however small $\Gamma$ is: those of Otto and Reznikoff \cite{OttoReznikoff07} and Marton \cite{Marton13} require pointwise bounds on the mixed second derivatives of the interaction, which are infinite, and the Zegarlinski coefficients \cite{Zegarlinski92} remain of order $1$ instead of $\frac1N$, so that Zegarlinski's condition fails for $N$ large. We explain this obstruction in Section~\ref{sec:obstruction}. The second mechanism, using some ideas from Rosenzweig \cite{Rosenzweig_partI}, avoids interdependence coefficients for the entropy. The relative entropy with respect to the free product measure (namely the Haar measure $\mathbf m_N$ on the torus)
  satisfies the independent tensorization, and the interaction only appears through an energy term, controlled by a uniform bound on the partition function, valid at every temperature, due to Wang and Zhao \cite{WangZhao26} (see also Duerinckx and Jabin \cite{DJ26} at high temperature). This produces a log-Sobolev inequality with an additive defect, uniform in $N$, which is removed using  Rothaus' lemma and the uniform Poincar\'e inequality obtained by the first mechanism.
\smallskip

The paper thus clarifies the scope of interdependence coefficients for singular interactions, identifies the obstacle to their use for the LSI, and determines under which assumptions on the interaction each propagation principle applies, see Theorems~\ref{T-abstract-W} and~\ref{T-abstract-LSI}.
We refer also to Section~\ref{sec.discu} compares the two propagation mechanisms and the role of each assumption.

\subsection{Main results}\label{subsec:main}

Our first main result concerns the uniform (in $N$) Poincar\'e and $W_2I$ inequalities. All the constants are explicit.

\bthm\label{T_W2I} Let $\mu_N$ be the Coulomb gas \eqref{Coulomb-measure} on $(\T^3)^N$, assume that $\Gamma<2\pi$, and recall that $W_*=\min_{\T^3}W$ $($numerically, $W_*=W(\frac12,\frac12,\frac12)\approx-0.802$, see Section \ref{sec.N}$)$.
\begin{enumerate}
  \item[(1)] If
  \bbeq\label{T_W2I1}
  r_{\rm D}:= \Gamma e^{-\Gamma W_*} \|W\|_{L^1(\T^3,dx)} <1,
  \nneq
  then the following uniform Poincar\'e inequality holds:
  \bbeq\label{T_PI}
  \var_{\mu_N} (f) \le  \frac{1}{(1-r_{\rm D})\cdot 2\pi(2\pi-\Gamma)} \int_{(\T^3)^N} |\nabla f|^2\, d\mu_N, \qquad f\in H^1(\mu_N).
  \nneq
  Equivalently, the spectral gap $\lambda_1(\mu_N)$ of $-\LL$ on $L^2(\mu_N)$ satisfies $\lambda_1(\mu_N) \geqslant (1-r_{\rm D})\cdot 2\pi(2\pi-\Gamma)$. As a consequence, the semigroup $P_t$ of the Langevin diffusion \eqref{Langevin_D} satisfies $\var_{\mu_N} (P_t f) \le  e^{-2 \lambda_1(\mu_N) t} \var_{\mu_N} (f)$ for all $f\in L^2(\mu_N)$ and $t\geqslant 0$.

  \item[(2)] If
  \bbeq\label{T_W2I3}
  r_{W_2}:=  \Gamma  e^{-\Gamma W_*/2} \sqrt{ \frac{7}{3}
       \left[1+ \frac{1+2\pi\Gamma}{2\pi(2\pi-\Gamma)} \right]} \cdot \sup_{e\in \rr^3,\, |e|=1}\| e\cdot \nabla W\|_{L^{6/5}(\T^3, dx)}<1,
  \nneq
  then the following uniform $W_2I$ inequality holds:
  \bbeq\label{T_W2I2}
  W_{2}^2(\nu, \mu_N) \le   [2C_{W_2I}]^2\,  I(\nu|\mu_N), \qquad \forall\, \nu\in \mathcal P_1 ((\T^3)^N),
  \nneq
  where
  \bbeq\label{T_W2I5}
  C_{W_2I} = \frac{C_{W_2I,1}}{1-r_{W_2}},\qquad
  C_{W_2I,1}= \sqrt{\left(\frac{3}{8}+ \frac{24}{2\pi(2\pi-\Gamma)}\right) \frac{1}{2\pi(2\pi-\Gamma)}}.
  \nneq
\end{enumerate}
\nthm

The assumption $\Gamma<2\pi$ is in fact implied by \eqref{T_W2I1} and by \eqref{T_W2I3}. Indeed, $W_*<0$ and, by Fourier series,
$$
\|W\|_{L^1}\geqslant\Big|\int_{\T^3}W(x)\cos(2\pi x_1)\,dx\Big|=\frac1\pi,\qquad \|\partial_1W\|_{L^{6/5}}\geqslant\Big|\int_{\T^3}\partial_1W(x)\sin(2\pi x_1)\,dx\Big|=2.
$$
Thus, $r_{\rm D}<1$ forces $\Gamma<\pi$ and $r_{W_2}<1$ forces $\Gamma<\frac{1}{2}\sqrt{3/7}<0.33$. We keep it in the statement since it is the condition under which the one-particle Poincar\'e inequality of Proposition~\ref{PI_single} holds. With the numerical values of Section \ref{sec.N}, \eqref{T_W2I1} holds for $\Gamma\le0.84$ and \eqref{T_W2I3} holds for $\Gamma\le0.13$.

The $W_2I$ inequality implies Talagrand's inequality, hence a Gaussian concentration inequality for $\mu_N$ which is dimension-free, and a lower bound on the ground state energy of Schr\"odinger operators $-\LL+v$, see Corollaries~\ref{cor1} and~\ref{cor2}. Theorem~\ref{T_W2I} extends to the two-dimensional Coulomb gas on $\T^2$ and to periodic Riesz gases $|x|^{-s}$, $0<s\le1$, on $\T^3$, see Theorems~\ref{T_2D} and~\ref{T_Riesz}.
With the same approach, the Poincaré inequality with a constant uniform in $N$ extends directly to all $d \ge 4$ (see the Section \ref{se.dim}). 

Our second main result is a uniform (in $N$) log-Sobolev inequality.

\bthm\label{T_LSI_intro}
Assume that $r_{\rm D}<1$, i.e.\ \eqref{T_W2I1} holds. Then, for every $N\geqslant2$ and every $f\in H^1(\mu_N)$,
$$
\Ent_{\mu_N}(f^2)\le 2C_{\rm LS}\int_{(\T^3)^N}|\nabla f|^2\,d\mu_N,
$$
where $C_{\rm LS}<\infty$ depends only on $\Gamma$ and is given in Theorem~\ref{T_LSI}. In particular, $C_{\rm LS}$ does not depend on $N$.
\nthm

%

The constant $C_{\rm LS}$ is explicit, up to the constant $\sup_{N\geqslant2}\log Z_N(3)$ of the partition bound, which is finite for every $\Gamma>0$ by the result of Wang and Zhao \cite{WangZhao26} (see Lemma~\ref{le.partition}). In particular, the only restriction on $\Gamma$ in Theorem~\ref{T_LSI_intro} is the condition $r_{\rm D}<1$ of the uniform Poincar\'e inequality. Moreover, by Theorem \ref{T_LSI_intro}, the entropy of the law of the Langevin dynamics \eqref{Langevin_D} with respect to $\mu_N$ decays exponentially fast, at the rate $2/C_{\rm LS}$, uniformly in $N$. 
 By the hierarchy recalled in Section~\ref{sec:prelim}, Theorem~\ref{T_LSI_intro} also yields uniform $W_2H$ and $W_2I$ inequalities. The two approaches we use in this work are complementary. The log-Sobolev inequality is the strongest of these functional inequalities, and it implies the others. The dependent tensorization of Theorem~\ref{T_W2I}\emph{(2)} reaches the $W_2I$ inequality directly, with simple explicit constants and without any bound on the partition function. 
Theorem~\ref{T_LSI_intro} also extends to the two-dimensional Coulomb gas and to periodic Riesz gases with $0<s\le1$, see Corollary~\ref{cor:LSI_ext}.

\subsection{Higher dimensions}
\label{se.dim}
Dimension $d=3$ is both the physical dimension and the largest dimension in which the Coulomb kernel is Hilbert-Schmidt, i.e.\ $W\in L^2(\T^d,dx)$, equivalently $\nabla W\in H^{-1}(\T^d,dx)$. Let $W$ be defined on $\T^d$ by $\hat W(0)=0$ and $\hat W(k)=\frac{1}{\pi|k|^2}$, $k\neq0$. For $d\geqslant4$, the situation is as follows.
\begin{itemize}
\item[-] The uniform Poincar\'e inequality of Theorem~\ref{T_W2I}\emph{(1)} holds on $\T^d$ for every $d\geqslant3$, with the same proof: the ground state transformation only uses $\Delta W=4\pi$ away from $0$, and the Dobrushin coefficients only use $W_*>-\infty$ and $W\in L^1(\T^d,dx)$.
\item[-] For every $d\geqslant2$, $\mu_N$ satisfies a $W_1H$ inequality with respect to the $\ell^1$-metric, with a constant of order $N$, and thus a Gaussian concentration inequality uniform in $N$, as soon as $\Gamma e^{-\Gamma W_*}\|\nabla W\|_{L^1(\T^d,dx)}$ is small enough (Proposition~\ref{P_W1H} in the appendix). 
\item[-] In the proof of Proposition~\ref{P-coef_W2}, the $W_2$-interdependence coefficients are bounded by $\frac{\Gamma}{N-1}$ times the $H^{-1}(\mu_t)$-norm of $e\cdot\nabla W(\cdot-\eta(t))$, and Lemma~\ref{lem26} bounds this norm by $\|e\cdot\nabla W\|_{L^{2d/(d+2)}}$ through the Sobolev embedding. For $d\geqslant4$, $\nabla W\notin L^{2d/(d+2)}(\T^d,dx)$, and this bound is infinite. Consequently, our approach does not provide a $W_2I$ inequality uniform in $N$ when $d\geqslant4$.
\item[-] By \cite[Proposition~2.2(ii)]{DJ26}, $Z_N(t)\to\infty$ for every $t>0$ when $W\notin L^2$ (which is the case when $d\ge 4$). In particular $H(\mathbf m_N|\mu_N)=\log Z_N(1)\to\infty$, and any argument comparing $\mu_N$ with the free product measure, such as Theorem~\ref{T-abstract-LSI}, fails.
\end{itemize}
The two mechanisms of this paper for the $W_2I$ and the log-Sobolev inequalities thus stop exactly at the same threshold, and the case $d=3$ is critical for them (see also the discussion at the end of Section~\ref{sec.discu}).


\subsection{Related literature}\label{subsec:literature}

 Functional inequalities with constants independent of the number of particles have been studied intensively, both for spin systems and for mean-field diffusions. For spin systems, the classical approach relies on Dobrushin-type weak-dependence conditions \cite{Dobrushin68,Dobrushin70}: Zegarlinski \cite{Zegarlinski92} derived the LSI, and the last named author \cite{Wu_AOP06} Poincar\'e and transport inequalities, under the Dobrushin uniqueness condition; Otto and Reznikoff \cite{OttoReznikoff07} and Marton \cite{Marton13} (see also \cite{CMT15}) gave criteria for the LSI based on bounds on the mixed second derivatives of the Hamiltonian. For mean-field systems, see \cite{Malrieu01,CMV03} in the convex case, \cite{GLWZ2022} in a perturbative regime (using the criterion of Ledoux \cite{Ledoux01b} for the Poincar\'e inequality and Zegarlinski's interdependence coefficients for the LSI), and \cite{DGPS23,BauerschmidtBodineau19,BBD25,CNZ24,Wang25} for more recent developments.

Much less is known for singular interactions. For Coulomb gases on $\rr^d$, $d\geqslant3$, with quadratic confinement, Becker and Menegaki \cite{Becker-Menegaki26} proved a uniform Poincar\'e inequality at weak coupling, through a ground state transformation and a compactness argument for the one-particle conditional laws, and Dobrushin's coefficients in total variation; their constants are not explicit. Theorem~\ref{T_W2I}\emph{(1)} is a quantitative counterpart of their result on the torus.  
 Rosenzweig \cite{Rosenzweig26} proved a uniform LSI for the two-dimensional Coulomb gas on $\rr^2$ at every temperature, after a perturbative proof \cite{Rosenzweig_partI} combining one-site LSI, a conditional entropy inequality, a partition bound \cite{DG25} and Rothaus' lemma. Our Section~\ref{sec:LSI} follows this perturbative strategy on $\T^3$, the key input being a uniform bound on the partition function for square-integrable interactions. Such a bound was obtained by Duerinckx and Jabin \cite[Theorem~2.1(i)]{DJ26} at high temperature, and very recently by Wang and Zhao \cite{WangZhao26} at every positive temperature, for logarithmic, Riesz and Bessel-Riesz interactions in the range $0\le s<\frac d2$, which includes the Coulomb kernel in dimension $3$. The latter result allows us to remove any smallness condition on $\Gamma$ coming from the partition function.
For background on Coulomb and Riesz gases, we refer to the monographs \cite{Chafai21,Serfaty26}.

\begin{remark}[Use of AI tools]\label{rem:AI}
AI tools were used only to help obtain the estimate of the weighted $H^{-1}(\mu^{\y})$-norm of $\nabla W(\cdot-y^j)$ in Lemma~\ref{lem26}, to compute the numerical constants of Section~\ref{sec.N}, and to understand why Zegarlinski's condition fails for the Coulomb gas (see more precisely the example given in the part dedicated to the Zegarlinski's coefficients at the end of  Section \ref{sec:obstruction}).   We have checked all these arguments and computations. The overall strategy, the other results and their proofs are due to the authors.
\end{remark}

\subsection{Organization of the paper}
Section~\ref{sec:W2I_all} is devoted to the dependent tensorization. After setting up the analytic framework for singular potentials (Section~\ref{subsec:FI}), we prove uniform Poincar\'e and $W_2I$ inequalities for one particle (Section~\ref{subsec:PI1}) and the key $H^{-1}$ estimate (Section~\ref{subsec:H-1}). We then tensorize the Poincar\'e inequality through Dobrushin's coefficients in total variation (Section~\ref{sec:tv}), and the $W_2I$ inequality through $W_2$-interdependence coefficients, which we estimate by the Hamilton-Jacobi equation (Section~\ref{sec:W2coef}). Section~\ref{sec:ext2} extends these results to the two-dimensional Coulomb gas and to periodic Riesz gases. Section~\ref{sec:LSI} is devoted to the LSI. We first explain why the previous mechanism does not apply to the entropy (Section~\ref{sec:obstruction}). We then prove a uniform one-site LSI (Section~\ref{sec:LSI1}), establish an abstract propagation principle with a defect (Section~\ref{sec:defect}), and combine it with the partition bound of \cite{DJ26} to prove Theorem~\ref{T_LSI_intro} (Sections~\ref{sec:partition} and~\ref{sec:LSI_proof}).  
Section~\ref{sec:LSI_ext} discusses extensions, and Section~\ref{sec.discu} compares the two propagation mechanisms and the role of each assumption.



\section{Dependent tensorization: uniform Poincar\'e and $W_2I$ inequalities}\label{sec:W2I_all}

In this section we prove Theorem~\ref{T_W2I}. We first set up the analytic framework for singular potentials (Section~\ref{subsec:FI}). We then work with a single particle on $\T^3$ in the field of $N-1$ fixed poles, namely with the conditional distributions \eqref{cond_law} (Sections~\ref{subsec:PI1} and~\ref{subsec:H-1}). The propagation to $\mu_N$ is the object of Sections~\ref{sec:abstractW}, \ref{sec:tv} and~\ref{sec:W2coef}, and Section~\ref{sec:ext2} is devoted to extensions.

\subsection{Langevin dynamics with singular potential and Fisher-Donsker-Varadhan information}\label{subsec:FI}
The potentials we consider are singular: they blow up at the poles, or at the collision set. In this subsection we show that this singularity is harmless. The Langevin dynamics is well posed from every point and never reaches the singularities (Proposition~\ref{prop:wellposed}). Moreover, its generator is determined by its action on smooth functions supported away from the singularities (Proposition~\ref{prop:Wu}). As a consequence, such functions are dense in the energy space $H^1(\mu)$ (Proposition~\ref{prop:Wu}\emph{(c)}). This density justifies the integrations by parts used throughout the paper, and gives a precise meaning to the Fisher-Donsker-Varadhan information \eqref{FI} (Definition~\ref{def:FI}).

Throughout, $V$ denotes either the one-particle potential on $\T^3$
\begin{equation}\label{eq.V1}
V(x)=\sum_{j=2}^{N}\gamma_jW(x-y^j), \qquad  \gamma_j>0,
\end{equation}
for a fixed configuration of poles $\y=(y^2,\cdots,y^N)\in(\T^3)^{N-1}$, or the $N$-particle potential on $(\T^3)^N\simeq\T^{3N}$
\begin{equation}\label{eq.V2}
V(x)=\sum_{i,j=1,i\neq j}^N\gamma_{ij}W(x^i-x^j), \qquad \gamma_{ij}>0.
\end{equation}
We write $\T^n$, with $n=3$ or $n=3N$, for the state space, $\mu=\frac1Ze^{-V}dx$, $\SS:=\{V=+\infty\}$ for the singular set (the set of poles $\{y^2,\dots,y^N\}$ in the first case, and the collision set $\SS_N$ in the second one), and $\mathcal D:=\T^n\setminus\SS$, which is open and connected. Since $W(x)=|x|^{-1}+R(x)$ near $0$ with $R$ smooth, and $W$ is smooth on $\T^3\setminus\{0\}$, the potential $V$ has the following properties.
\begin{enumerate}[(\bfseries P1)]
  \item $e^{-V}$ is smooth on $\T^n$ (near $\SS$ it behaves like $e^{-\gamma/|z|}$, all of whose derivatives vanish at $z=0$). In particular the ground state $\varphi:=Z^{-1/2}e^{-V/2}$ belongs to $H^1(\T^n,dx)$ and $\mu=\varphi^2dx$.
  \item $V\in C^\infty(\mathcal D)$, $V$ is bounded from below since $W\geqslant W_*$, and $V(x)\to+\infty$ as $x\to\SS$ (one term blows up while the others remain bounded from below). Consequently, the sublevel sets $\{V\le  k\}$, $k\in\nn$, are compact subsets of $\mathcal D$, on which $\nabla V$ is Lipschitz.
  \item $\Delta V$ is constant on $\mathcal D$: since $\Delta W=4\pi$ on $\T^3\setminus\{0\}$ by \eqref{eq:W},
  $$
  \Delta V=4\pi\sum_{j=2}^N\gamma_j\quad\text{for \eqref{eq.V1}},\qquad
  \Delta V=8\pi\sum_{i\neq j}\gamma_{ij}\quad\text{for \eqref{eq.V2}},
  $$
  where in the second case each term $W(x^i-x^j)$ is differentiated both in $x^i$ and in $x^j$.
\end{enumerate}

\subsubsection{Well-posedness of the Langevin dynamics}
Consider the Langevin diffusion on $\T^n$
\bbeq\label{SDE_V}
dX_t=-\nabla V(X_t)\,dt+\sqrt2\,dB_t,
\nneq
where $B_t$ is a standard Brownian motion in $\rr^n$, and its generator
$\LL f=\Delta f-\nabla V\cdot\nabla f$, $f\in C_c^\infty(\mathcal D)$.
A weak solution of \eqref{SDE_V} in $\mathcal D$ is a continuous $\mathcal D$-valued process satisfying \eqref{SDE_V} for some Brownian motion on some filtered probability space.

\bprop\label{prop:wellposed}
For every initial point $x\in \mathcal D$:
\begin{enumerate}[(a)]
  \item The SDE \eqref{SDE_V} has a unique strong solution $(X_t^x)_{t\geqslant0}$, which is defined for all times and a.s.\ never reaches $\SS$.
  \item Uniqueness in law holds for weak solutions of \eqref{SDE_V} in $\mathcal D$. Thus, the martingale problem associated with $(\LL,C_c^\infty(\mathcal D))$ is well posed on $\mathcal D$.
\end{enumerate}
Consequently, the laws $(\mathbb P_x)_{x\in \mathcal D}$ of these solutions form a conservative strong Markov process on $\mathcal D$, whose transition semigroup we denote by $P_tf(x):=\mathbb E_x f(X_t)$.
\nprop

\bprf
(a) By (\textbf{P2}), \eqref{SDE_V} has a unique strong solution up to the exit time $\tau_k$ of $\{V\le  k\}$, for each $k$, and thus up to $\tau:=\lim_k\tau_k$, the hitting time of $\SS$. By (\textbf{P3}), $\LL V=\Delta V-|\nabla V|^2\le  c$ on $\mathcal D$ for a constant $c$, and $V$ is a Lyapunov function. It\^o's formula gives $\mathbb E\,V(X_{t\wedge\tau_k})\le  V(x)+ct$, and since $V\geqslant\inf V>-\infty$, $\mathbb P(\tau_k\le  t)\le (V(x)-\inf V+ct)/(k-\inf V)\to0$ as $k\to\infty$. Therefore $\tau=+\infty$ a.s. This is Khasminskii's argument.

(b) By (a), pathwise uniqueness holds, and uniqueness in law follows from the Yamada-Watanabe theorem, see \cite[Lemma~18.17]{Kallenberg97}. The well-posed martingale problem then defines a conservative strong Markov family, see \cite[Theorem~18.11]{Kallenberg97}.
\nprf

In particular, for the Coulomb gas $\mu_N$, where $\gamma_{ij}=\frac{\Gamma}{2(N-1)}$ and $\Delta V=4\pi\Gamma N$, no collision occurs in the dynamics \eqref{Langevin_D}.

\subsubsection{Link with the Dirichlet form}
By integration by parts, $\LL$ is symmetric on $L^2(\mu)$:
$$
\<-\LL f,g\>_\mu=\<f,-\LL g\>_\mu=\int_{\T^n}\nabla f\cdot\nabla g\,d\mu,\qquad f,g\in C_c^\infty(\mathcal D).
$$
Two Dirichlet forms on $L^2(\mu)$ extend $(\EE_\mu,C_c^\infty(\mathcal D))$, where $\EE_\mu[g]=\int_{\mathcal D}|\nabla g|^2d\mu$. The \emph{minimal} one, $(\EE_\mu,\mathcal F_{\min})$, is the closure of $(\EE_\mu,C_c^\infty(\mathcal D))$, which is closable \cite{AHS77,FOT94,Albeverio03}. The \emph{maximal} one is $(\EE_\mu,H^1(\mu))$, where
\bbeq\label{H1mu}
H^1(\mu)=\left\{g\in L^2(\mu);\ \nabla g\in L^1_{loc}(\mathcal D,dx),\ \EE_\mu[g]:=\int_{\mathcal  D}|\nabla g|^2\,d\mu<+\infty\right\},
\nneq
and $\nabla g$ is taken in the distributional sense on $\mathcal D$. It is closed because the density of $\mu$ is bounded from below on every compact subset of $\mathcal D$, and Markovian because normal contractions decrease $|\nabla g|$. The next proposition shows that the process of Proposition~\ref{prop:wellposed} is the diffusion associated with these forms, and that the two forms coincide.

\bprop\label{prop:Wu}The following hold:
\smallskip
\begin{enumerate}[(a)]
  \item {\rm(Identification)} The transition semigroup $(P_t)$ of Proposition~\ref{prop:wellposed} is the $L^2(\mu)$-semigroup $(e^{t\LL_{\min}})$ associated with $(\EE_\mu,\mathcal F_{\min})$: for every bounded measurable $f$ and $t\geqslant0$, $P_tf=e^{t\LL_{\min}}f$ $\mu$-a.e. In particular $\mu$ is reversible for $(X_t)$.
  \item {\rm($L^1$-uniqueness)} $(\LL,C_c^\infty(\mathcal D))$ is an essential generator in $L^1(\mu)$: $(P_t)$ is the unique strongly continuous semigroup on $L^1(\mu)$ whose generator extends $(\LL,C_c^\infty(\mathcal D))$.
  \item {\rm(Uniqueness of the Dirichlet form)} $(\EE_\mu,C_c^\infty(\mathcal D))$ has a unique Dirichlet form extension. In particular $\mathcal F_{\min}=H^1(\mu)$, i.e.\ $C_c^\infty(\mathcal D)$ is dense in $H^1(\mu)$: for every $g\in H^1(\mu)$ there exists a sequence $(g_k)$ in $C_c^\infty(\mathcal D)$ such that $\int_{\mathcal D}\big((g_k-g)^2+|\nabla g_k-\nabla g|^2\big)\,d\mu\to0$.
\end{enumerate}
\nprop

\bprf
We apply \cite{Wu_PTRF99} with $\varphi=Z^{-1/2}e^{-V/2}$, which is continuous on $\T^n$, strictly positive on $\mathcal D$, with $\nabla\varphi\in L^2(\T^n,dx)$ by (\textbf{P1}), and with $|\nabla\varphi|^2/\varphi^2=\frac14|\nabla V|^2$ locally bounded on $\mathcal D$ by (\textbf{P2}). Thus, the assumptions (\textbf{H1}), (\textbf{H2}) and (0.2) of \cite{Wu_PTRF99} hold.\footnote{The results of \cite{Wu_PTRF99} are stated on domains of $\rr^d$, with generator $\frac12\Delta+({\nabla\varphi}/{\varphi})\cdot\nabla= \LL/2$. They extend to $\T^n$, and the factor $\frac12$ only amounts to the time change $t\mapsto 2t$.}

\emph{(a)} By \cite[Lemma~1.7]{Wu_PTRF99}, for every initial law $\nu\ll\mu$, the weak solution of \eqref{SDE_V} killed at $\partial \mathcal D=\SS$ and satisfying the local finite energy condition $\int_0^t|\nabla V|^2(X_s)\,ds<\infty$ for $t<\tau$ is unique, and it is the Markov process associated with $(\EE_\mu,\mathcal F_{\min})$. The solution of Proposition~\ref{prop:wellposed} satisfies these requirements: it never reaches $\SS$, and its paths are continuous in $\mathcal D$, which implies that $|\nabla V|^2(X_s)$ is bounded on $[0,t]$.

\emph{(b)} By \emph{(a)} and Proposition~\ref{prop:wellposed}\emph{(a)}, the diffusion associated with $(\EE_\mu,\mathcal F_{\min})$ a.s.\ never reaches $\SS$, which is condition (2.1) of \cite{Wu_PTRF99}. By \cite[Lemma~1.2]{Wu_PTRF99} and \emph{(a)}, $(P_t)$ is strongly continuous on $L^1(\mu)$ and its generator extends $(\LL,C_c^\infty(\mathcal D))$. Therefore \cite[Theorem~2.1 and Proposition~1.1]{Wu_PTRF99} show that the closure of $(\LL,C_c^\infty(\mathcal D))$ in $L^1(\mu)$ is the generator of $(P_t)$ in $L^1(\mu)$, and that $(P_t)$ is the unique strongly continuous semigroup on $L^1(\mu)$ whose generator extends $(\LL,C_c^\infty(\mathcal D))$. See also \cite[Corollary~2.2(a)]{Wu_PTRF99}.

\emph{(c)} Let $(\mathcal E,\mathcal F)$ be a Dirichlet form on $L^2(\mu)$ whose $L^2(\mu)$-generator $A$ extends $(\LL,C_c^\infty(\mathcal D))$. Its semigroup $(e^{tA})$ is symmetric and sub-Markovian, and by duality $\|e^{tA}f\|_{L^1(\mu)}\le \|f\|_{L^1(\mu)}$. Since $\mu$ is finite, $L^2(\mu)$ is dense in $L^1(\mu)$, and $(e^{tA})$ extends to a strongly continuous semigroup on $L^1(\mu)$. Its $L^1(\mu)$-generator also extends $(\LL,C_c^\infty(\mathcal D))$, since for $f\in C_c^\infty(\mathcal D)$, $t^{-1}(e^{tA}f-f)\to\LL f$ in $L^2(\mu)$, and thus in $L^1(\mu)$. By \emph{(b)}, this extension is $(P_t)$. As a result, $(e^{tA})$, and thus $A$ and $(\mathcal E,\mathcal F)$, are uniquely determined, see \cite[Corollary~2.2(a)]{Wu_PTRF99}. Both $(\EE_\mu,\mathcal F_{\min})$ and $(\EE_\mu,H^1(\mu))$ are such Dirichlet forms: for $f\in C_c^\infty(\mathcal D)$ and $g$ in the domain of the Dirichlet form, integration by parts on $\mathcal D$ (in the distributional sense) gives $\EE_\mu(f,g)=-\int_{\mathcal D}(\LL f)\,g\,d\mu$. Consequently, they coincide, i.e.\ $\mathcal F_{\min}=H^1(\mu)$.
\nprf

By \cite[Theorem~3.1 and Proposition~3.4]{Wu_PTRF99}, the identity in \emph{(a)} holds for every $x\in \mathcal D$, and $(P_t)$ is strong Feller on $\mathcal D$. Part~\emph{(c)} is the density result we shall use throughout to extend integrations by parts from $C_c^\infty(\mathcal D)$ to $H^1(\mu)$.

\bdef\label{def:FI}
The Fisher-Donsker-Varadhan information of $\nu\in\mathcal P_1(\T^n)$ with respect to $\mu$ is
$I(\nu|\mu)=\EE_\mu[h]$ if $\nu=h^2\mu$ with $h=\sqrt{d\nu/d\mu}\in H^1(\mu)$, and $I(\nu|\mu)=+\infty$ otherwise.
\ndef
For $\mu=\mu_N$ this is exactly \eqref{FI}, and by Proposition~\ref{prop:Wu}\emph{(c)} it agrees with the definition recalled in Section~\ref{sec:prelim}.

\subsection{Poincar\'e and $W_2I$ inequalities for one particle}\label{subsec:PI1}

Fix poles $\y=(y^2,\cdots,y^N)\in(\T^3)^{N-1}$ with charges $\gamma_j>0$, let $V^{\y}$ be the one-particle potential \eqref{eq.V1}, with singular set $\SS=\{y^2,\cdots,y^N\}$ and $\mathcal D=\T^3\setminus\SS$, and consider the one-particle distribution
\bbeq\label{mu_single}
\mu^{\y}(dx) = \frac{1}{Z(\y)}\, e^{-V^{\y}(x)}\, dx .
\nneq
This is the measure $\mu$ of Section~\ref{subsec:FI} for the potential \eqref{eq.V1}, and all results of that subsection apply to it. When $\gamma_j=\frac{\Gamma}{N-1}$ for all $j$ and $\y=x^{-i}$, the measure $\mu^{\y}$ is the conditional distribution $\mu_i(\cdot\,|\,x^{-i})$ of \eqref{cond_law}. We denote by $\bar\gamma:=\sum_{j=2}^{N}\gamma_j$ the total charge of the poles (for the Coulomb gas, $\bar\gamma=\Gamma$), and by $\LL^{\y}=\Delta-\nabla V^{\y}\cdot\nabla$ the Langevin generator associated with $\mu^{\y}$. Note that, since $\int_{\T^3}V^{\y}dx=0$ and $V^{\y}\geqslant\bar\gamma W_*$, Jensen's inequality gives
\bbeq\label{eq:Zbound}
1\le Z(\y)\le e^{-\bar\gamma W_*},\qquad\text{and thus}\qquad \frac{d\mu^{\y}}{dx}\le e^{-\bar\gamma W_*}.
\nneq

\bprop\label{PI_single} If the total charge satisfies $\bar\gamma < 2\pi$, then
\bbeq\label{PI_single2}
\var_{\mu^{\y}}(f) \le  C_{{\rm P},1} \int_{\T^3} |\nabla f|^2\, d\mu^{\y}, \qquad f\in H^1(\mu^{\y}),\qquad C_{{\rm P},1}=\frac{1}{ 2\pi (2\pi -\bar\gamma)}.
\nneq
The constant $C_{{\rm P},1}$ depends only on $\bar\gamma$; it is independent of $N$, of the number of poles and of their positions. Equivalently, the spectral gap of $\LL^{\y}$ on $L^2(\mu^{\y})$ is at least $2\pi(2\pi -\bar\gamma)$.
\nprop

This is the torus counterpart of \cite[Theorem~2.7]{Becker-Menegaki26}, which is proved in $\rr^d$ with a Gaussian confinement by a ground state transformation and a compactness argument. On $\T^3$, the ground state transformation combined with the Poincar\'e inequality for the Lebesgue measure avoids the compactness argument and gives an explicit constant.

\bprf Write $\mu=\mu^{\y}$ and $V=V^{\y}$.

{\bf Step 1. Ground state transformation.}
Since $\Delta W=4\pi$ on $\T^3\setminus\{0\}$ by \eqref{eq:W},
\begin{equation}\label{eq:Delta}
 \Delta V(x) = \sum_{j=2}^{N} \gamma_j\, \Delta W(x-y^j)  =4\pi\bar\gamma,\qquad x\in \mathcal D.
\end{equation}
Let $\varphi=Z^{-1/2}e^{-V/2}$ be the ground state, so that $\mu=\varphi^2dx$. Since $\int f^2\,d\mu=\int(f\varphi)^2\,dx$, the ground state transformation $f\mapsto f\varphi$ is an isometry from $L^2(\mu)$ onto $L^2(\T^3,dx)$. Let $f\in C^\infty_c(\mathcal D)$ and $g:=f\varphi\in C^\infty_c(\mathcal D)$. Since $\nabla\varphi=-\frac12\varphi\nabla V$ on $\mathcal D$, we have $\varphi\nabla f=\nabla g+\frac12 g\nabla V$. Expanding the square, and integrating by parts in the last term (which is justified because $g^2$ has compact support in $\mathcal D$, where $V$ is smooth), we get by \eqref{eq:Delta}:
\begin{align*}
 \EE_\mu[f]=\int_{\T^3}|\nabla f|^2\varphi^2\,dx
 &=\int_{\T^3} \left|\nabla g+\tfrac12g\nabla V\right|^2 dx \\
 &=\int_{\T^3} |\nabla g|^2 dx +\tfrac14\int_{\T^3}|\nabla V|^2g^2\, dx
    +\tfrac12\int_{\T^3}\nabla(g^2)\cdot\nabla V\, dx \\
 &=\int_{\T^3} |\nabla g|^2 dx  +\tfrac14\int_{\T^3} |\nabla V|^2g^2\, dx -\tfrac12\int_{\T^3}g^2\,\Delta V\,dx\\
 &=\int_{\T^3} |\nabla g|^2 dx  +\tfrac14\int_{\T^3} |\nabla V|^2g^2\, dx -2\pi\bar\gamma \int_{\T^3} g^2\, dx .
\end{align*}
Dropping the nonnegative term $\frac14\int|\nabla V|^2g^2\,dx$, we get that:
\begin{equation}\label{eq:H1}
      \int_{\T^3} |\nabla (f\varphi )|^2\, dx \le  \EE_\mu[f]+2\pi\bar\gamma \int_{\T^3} f^2\, d\mu.
\end{equation}
We now extend \eqref{eq:H1} to all $f\in H^1(\mu)$. By Proposition~\ref{prop:Wu}\emph{(c)}, there are $f_k\in C_c^\infty(\mathcal D)$ with $f_k\to f$ in $H^1(\mu)$. Then, $f_k\varphi\to f\varphi$ in $L^2(\T^3,dx)$, and \eqref{eq:H1} applied to $f_k-f_l$ shows that $(\nabla(f_k\varphi))$ is Cauchy in $L^2(\T^3,dx)$. Therefore $f\varphi\in H^1(\T^3,dx)$, and \eqref{eq:H1} holds for $f$ by passing to the limit.

{\bf Step 2. Poincar\'e inequality.}
Let $f\in H^1(\mu)$ and set $\tilde{f} = f- c_f$, where $c_f:=\int f\varphi\,dx/\int\varphi\,dx$ (note that $\int\varphi\,dx>0$). Set also $\tilde g := \tilde{f} \varphi$. By the choice of $c_f$, $\tilde{g}$ is centered with respect to $dx$. Because the variance minimizes $a\mapsto\int(f-a)^2d\mu$,
\bbeq\label{3eq7}
\var_{dx}(\tilde{g})=\int_{\T^3} \tilde{g}^2\, dx = \int_{\T^3} (f-c_f)^2\, d\mu\geqslant \var_\mu(f).
\nneq
The Lebesgue measure on $\T^3$ satisfies the Poincar\'e inequality $\var_{dx}(g)\le \frac1{4\pi^2}\int|\nabla g|^2dx$, since $4\pi^2$ is the smallest nonzero eigenvalue of $-\Delta$ on the unit torus. Applying it to $\tilde g$, and then \eqref{eq:H1} to $\tilde{f}$ (note that $\EE_\mu[\tilde f]=\EE_\mu[f]$ and $\int\tilde f^2d\mu=\int\tilde g^2dx$), we get
$$
4\pi^2 \int_{\T^3} \tilde{g}^2\, dx \le  \int_{\T^3} |\nabla \tilde{g}|^2\, dx \le  \EE_\mu[f]+2\pi\bar\gamma \int_{\T^3} \tilde{g}^2\, dx.
$$
Thus, if $\bar\gamma<2\pi$, we deduce by \eqref{3eq7} that:
$
  2\pi(2\pi -\bar\gamma) \var_\mu(f) \le  2\pi(2\pi -\bar\gamma) \int_{\T^3} \tilde{g}^2\, dx \le   \EE_\mu[f].
$
\nprf

To pass from the Poincar\'e inequality to a $W_2I$ inequality, we use the following centering lemma of Liu.

\blem\label{le.Liu} {\rm (\cite[Theorem~1.5 and Remark~1.6]{LiuY_EJP})} Let $\mu$ be a probability measure on $\T^n$, absolutely continuous with respect to $dx$, satisfying the Poincar\'e inequality $\var_\mu(g)\le  C_P\int_{\T^n}|\nabla g|^2\,d\mu$ for every Lipschitz function $g$. For $\nu=h^2 \mu\neq \mu$ with $h\geqslant 0$, let
$
\nu_c := \frac{(h-\mu(h))^2}{\var_{\mu}(h)}\, \mu
$
be the associated centered probability measure. Then, we have:
\bbeq\label{LiuY}
W_2^2(\nu, \mu) \le  2 \var_\mu(h)\, W_2^2(\nu_c, \mu) + 96C_P \var_\mu(h).
\nneq
\nlem
In \cite{LiuY_EJP}, this result is stated for $\mu=e^{-V}dx$ with $V\in C^1$ on a complete Riemannian manifold, which does not cover $\mu^{\y}$, since $V^{\y}=+\infty$ at the poles. However, its proof \cite[Section~6]{LiuY_EJP} uses neither the generator nor the regularity of $V$: it only uses the Poincar\'e inequality for the Lipschitz functions $Q_tv-\mu(Q_tv)$ and the properties of $Q_t$ collected in Lemma~\ref{le.HJ} below (in \cite{LiuY_EJP}, $\frac1t$ is replaced by $\frac1{2t}$, which only amounts to a time change). It thus applies to $\mu=\mu^{\y}$ with $C_P=C_{{\rm P},1}$, since Lipschitz functions belong to $H^1(\mu^{\y})$.

Lemma~\ref{le.Liu} yields a $W_2I$ inequality as follows: on the compact space $\T^n$, the transport cost $W_2^2(\nu_c,\mu)$ is bounded by the squared diameter, and the variance $\var_\mu(h)$ is controlled by the Fisher-Donsker-Varadhan information through the Poincar\'e inequality. More precisely, if $\mu$ satisfies the assumptions of Lemma~\ref{le.Liu} on $\T^n$, then
\bbeq\label{W2I_dim}
W_2^2(\nu, \mu) \le  \left(\frac{n}{2} + 96C_P\right) C_P\, I(\nu|\mu),\qquad \forall\,\nu\in\mathcal P_1(\T^n).
\nneq
Indeed, if $I(\nu|\mu)=+\infty$ or $\nu=\mu$ there is nothing to prove. Otherwise, $\nu=h^2\mu$ with $h\in H^1(\mu)$, $h\geqslant0$. Since every coordinate of $x-y$ can be taken in $[-\frac12,\frac12]$, we have $W_2^2(\nu_c, \mu)\le  {\rm Diam}(\T^n)^2=n/4$. Moreover $\var_\mu(h)\le  C_P\,\EE_\mu[h]=C_P\, I(\nu|\mu)$, and \eqref{LiuY} yields \eqref{W2I_dim}.
The constant in \eqref{W2I_dim} grows with the dimension. For $n=3N$ it cannot be applied directly to $\mu_N$, even though $\mu_N$ satisfies a uniform Poincar\'e inequality. It does apply, however, to the one-particle distribution $\mu^{\y}$ on $\T^3$. With $n=3$ and $C_P=C_{{\rm P},1}$, we have $\big(\frac32+96C_{{\rm P},1}\big)C_{{\rm P},1}=4\big(\frac38+24C_{{\rm P},1}\big)C_{{\rm P},1}$, which yields the following result.

\bprop\label{W2I_1particle} Assume that $\bar\gamma <2\pi$, and let
\bbeq\label{W2I_one particle1}
C_{W_2I,1}= \sqrt{\left(\frac{3}{8}+ 24 C_{{\rm P},1}\right) C_{{\rm P},1}} = \sqrt{\left(\frac{3}{8}+ \frac{24}{2\pi(2\pi-\bar\gamma)}\right) \frac{1}{2\pi(2\pi-\bar\gamma)}} .
\nneq
Then, for every configuration of poles $\y$, the one-particle distribution $\mu^{\y}$ satisfies
\bbeq\label{W2I_1particle2}
W_2^2(\nu, \mu^{\y}) \le  [2 C_{W_2I,1}]^2\, I(\nu|\mu^{\y}), \qquad \forall\, \nu\in \mathcal P_1(\T^3).
\nneq
\nprop

\subsection{$H^{-1}(\mu^{\y})$-norm estimate of $\nabla W(\cdot-y^j)$}\label{subsec:H-1}

In this subsection we prove the estimate which will be used in Section~\ref{sec:W2coef} to control the $W_2$-interdependence coefficients. It is the place where the singularity of the interaction is measured. We keep the notation of Section~\ref{subsec:PI1}. Fix $j\in\{2,\dots,N\}$ and consider the (repulsive) Coulomb force
$G^{{\scriptscriptstyle (j)}} (x):=-\nabla W(x-y^j)$
exerted by the pole $y^j$. Although $|G^{{\scriptscriptstyle (j)}}|^2$ is not integrable with respect to $dx$ near $y^j$ (since $|G^{{\scriptscriptstyle (j)}}(x)|\sim|x-y^j|^{-2}$), we have $G^{{\scriptscriptstyle (j)}}\in L^2(\mu^{\y})$. Indeed, $V^{\y}(x)\geqslant\gamma_jW(x-y^j)+(\bar\gamma-\gamma_j)W_*$ and $W(z)\geqslant|z|^{-1}-C$ near $0$, and the density of $\mu^{\y}$ is bounded by $C'e^{-\gamma_j/|x-y^j|}$ near $y^j$. As a result, $|G^{{\scriptscriptstyle (j)}}|^2$ times this density is bounded near $y^j$, while $G^{{\scriptscriptstyle (j)}}$ is bounded away from $y^j$.

By the Poincar\'e inequality \eqref{PI_single2}, $\EE_{\mu^{\y}}$ is coercive on the closed subspace $H^1_c(\mu^{\y}):=\{\phi\in H^1(\mu^{\y}),\ \mu^{\y}(\phi)=0\}$. Moreover, for $\phi\in H^1_c(\mu^{\y})$,
$$
\Big|\int_{\T^3}\big(G_k^{\scriptscriptstyle (j)}-\mu^{\y}(G_k^{\scriptscriptstyle (j)})\big)\,\phi\,d\mu^{\y}\Big|\le \|G_k^{\scriptscriptstyle (j)}\|_{L^2(\mu^{\y})}\,\sqrt{C_{{\rm P},1}\,\EE_{\mu^{\y}}[\phi]} .
$$
By the Riesz representation theorem, for $k=1,2,3$ there is a unique $u_k^{\scriptscriptstyle (j)}\in H^1_c(\mu^{\y})$ solving the Poisson equation
\begin{equation}\label{eq.poisson}
-\LL^{\y} u_k^{\scriptscriptstyle (j)} = G_k^{\scriptscriptstyle (j)} - \mu^{\y}(G_k^{\scriptscriptstyle (j)})
\end{equation}
in the weak sense, i.e.\
$
\int\nabla u_k^{\scriptscriptstyle (j)}\cdot\nabla\phi\,d\mu^{\y}=\int\big(G_k^{\scriptscriptstyle (j)}-\mu^{\y}(G_k^{\scriptscriptstyle (j)})\big)\,\phi\,d\mu^{\y}
$
for all $\phi\in H^1(\mu^{\y})$ (both sides vanish on constants). We write $\bu^{\scriptscriptstyle (j)}=(u_1^{\scriptscriptstyle (j)},u_2^{\scriptscriptstyle (j)},u_3^{\scriptscriptstyle (j)})$. The quantity $\int|\nabla(e\cdot\bu^{\scriptscriptstyle (j)})|^2d\mu^{\y}$ is the square of the $H^{-1}(\mu^{\y})$-norm of the centered force $e\cdot(G^{\scriptscriptstyle (j)}-\mu^{\y}(G^{\scriptscriptstyle (j)}))$.

\blem\label{lem26} Assume that $\bar\gamma<2\pi$. For every $j\in\{2,\dots,N\}$ and every unit vector $e\in \rr^3$,
\bbeq\label{lem26a}
\int_{\T^3} |\nabla (e\cdot \bu^{\scriptscriptstyle (j)} )|^2\,d\mu^{\y} \le  e^{-\bar\gamma W_*} C_{{\rm So}}
       \big[1+(1+2\pi\bar\gamma)C_{{\rm P},1}\big] \|e\cdot \nabla W\|^2_{L^{6/5}(\T^3,dx)},
\nneq
where $C_{{\rm P},1}=\frac{1}{2\pi(2\pi-\bar\gamma)}$ and $C_{{\rm So}}$ is the Sobolev constant in
\bbeq\label{SobolevI}
 \|w\|_{6}^2\le  C_{{\rm So}}
       \left(\|\nabla w\|_{2}^2+\|w\|_{2}^2\right),\qquad w\in H^1(\T^3,dx).
\nneq
In particular, the right-hand side of \eqref{lem26a} does not depend on $j$, $N$, nor on the positions of the poles.
\nlem

By Lemma \ref{le.Sobolev}, $C_{{\rm So}}\le 7/3$.  

\bprf Let $\varphi=Z^{-1/2}e^{-V^{\y}/2}$ be the ground state of $\mu^{\y}$, and set $f:=e\cdot\bu^{\scriptscriptstyle (j)} \in H^1_c(\mu^{\y})$. Applying \eqref{eq:H1} to $f$ and then the Poincar\'e inequality \eqref{PI_single2}, we get
$$
 \int_{\T^3}|\nabla(f\varphi)|^2\,dx
 \le  \int_{\T^3}|\nabla f|^2\,d\mu^{\y}
   +2\pi\bar\gamma \int_{\T^3} f^2\,d\mu^{\y}
 \le (1+2\pi\bar\gamma\, C_{{\rm P},1})
                     \int_{\T^3}|\nabla f|^2\,d\mu^{\y}.
$$
By the Sobolev inequality \eqref{SobolevI} and $\|f\varphi\|_2^2=\int f^2\,d\mu^{\y}\le  C_{{\rm P},1}\int|\nabla f|^2\,d\mu^{\y}$, we have:
\bbeq\label{2Eq17}
 \|f \varphi\|_6^2
 \le  C_{{\rm So}}\left(\|\nabla (f \varphi)\|_2^2+\|f \varphi\|_2^2\right)
 \le  C_{{\rm So}} \big[1+(1+2\pi\bar\gamma )C_{{\rm P},1}\big]
                         \int_{\T^3}|\nabla f|^2\,d\mu^{\y}.
\nneq
Multiplying the weak formulation of \eqref{eq.poisson} by $e_k$ and summing over $k$, we get
$
\int\nabla f\cdot\nabla\phi\,d\mu^{\y}=\int\big(e\cdot G^{{\scriptscriptstyle (j)}}-\mu^{\y}(e\cdot G^{{\scriptscriptstyle (j)}})\big)\,\phi\,d\mu^{\y}
$ for every $\phi\in H^1(\mu^{\y})$.
We take $\phi=f$. Since $\mu^{\y}(f)=0$, the constant $\mu^{\y}(e\cdot G^{{\scriptscriptstyle (j)}})$ does not contribute. By H\"older's inequality with exponents $\frac65$ and $6$, by $\varphi^2\le e^{-\bar\gamma W_*}$ (see \eqref{eq:Zbound}) and by \eqref{2Eq17}, we obtain
\begin{align*}
 \int_{\T^3}|\nabla f|^2\,d\mu^{\y}&=\int_{\T^3} \big((e\cdot  G^{{\scriptscriptstyle (j)}} )\varphi\big)\,(f\varphi)\,dx
 \le  \|(e\cdot  G^{{\scriptscriptstyle (j)}} ) \varphi\|_{6/5}\,\|f \varphi\|_6\\
& \le  e^{-\bar\gamma W_*/2}\|e\cdot  G^{{\scriptscriptstyle (j)}} \|_{6/5}
 \left\{C_{{\rm So}}\big[1+(1+2\pi\bar\gamma )C_{{\rm P},1}\big]
                    \int_{\T^3}|\nabla f|^2\,d\mu^{\y}\right\}^{1/2}.
\end{align*}
Dividing by the square root of $\int|\nabla f|^2\,d\mu^{\y}$ (if it is nonzero; otherwise there is nothing to prove), and using $\|e\cdot  G^{{\scriptscriptstyle (j)}} \|_{6/5}=\|e\cdot \nabla W\|_{6/5}$ (translation invariance of $dx$), we get \eqref{lem26a}.
\nprf

\subsection{Dependent tensorization: an abstract statement}\label{sec:abstractW}
Let $E=\T^d$ with its flat metric $d$, and let $\mu$ be a probability measure on $E^N$ with conditional distributions $\mu_i(\cdot|x^{-i})$, $1\le i\le N$. For a metric $\delta$ on $E$ (in this paper, $\delta$ is the discrete metric $\delta_0(x,y)=1_{x\neq y}$ or the flat metric $d$), and for a cost exponent $p\in\{1,2\}$, the interdependence coefficient of $x^i$ on $x^j$ is
$$
c_{ij}^{\delta,p}:=\sup_{x^{-i}=y^{-i}\text{ off }j}\frac{W_{p,\delta}\big(\mu_i(\cdot|x^{-i}),\mu_i(\cdot|y^{-i})\big)}{\delta(x^j,y^j)},\qquad c_{ii}^{\delta,p}:=0,
$$
where $W_{p,\delta}$ is the Wasserstein distance associated with $\delta$, and ``$x^{-i}=y^{-i}$ off $j$'' means $x^k=y^k$ for all $k\neq i,j$, and $x^j\neq y^j$. The coefficient $c_{ij}^{\delta,p}$ is the Lipschitz constant of $x^j\mapsto\mu_i(\cdot|x^{-i})$, uniformly in the other coordinates. The following statement summarizes the dependent tensorization principle. Part~(a) is due to the last named author \cite[Theorem~2.1, (2.4)]{Wu_AOP06}, and part~(b) to Gao and the last named author \cite[Theorem~3.1]{GW06}.

\bthm\label{T-abstract-W}
Let $\mu$ be a probability measure on $E^N$, with conditional distributions $\mu_i(\cdot|x^{-i})$.
\begin{enumerate}[(a)]
  \item {\rm (Poincar\'e inequality)} Assume that each $\mu_i(\cdot|x^{-i})$ satisfies a Poincar\'e inequality with a constant $C_{{\rm P},1}$ independent of $i$ and $x^{-i}$, and that the spectral radius $r$ of the matrix $(c^{\delta,1}_{ij})$ satisfies $r<1$ for some metric $\delta$. Then, $\mu$ satisfies a Poincar\'e inequality with constant $\frac{C_{{\rm P},1}}{1-r}$.
  \item {\rm ($W_2I$ inequality)} Assume that each $\mu_i(\cdot|x^{-i})$ satisfies a $W_2I$ inequality with a constant $C_{W_2I,1}$ independent of $i$ and $x^{-i}$, that the Fisher-Donsker-Varadhan information is additive in the sense of \eqref{eq:additivity} below, and that the operator norm $\rho$ of the matrix $C^{(2)}=(c^{d,2}_{ij})$ on $\ell^2$ satisfies $\rho<1$. Then, $\mu$ satisfies a $W_2I$ inequality with constant $\frac{C_{W_2I,1}}{1-\rho}$.
\end{enumerate}
\nthm

In both cases the constant of $\mu$ is the one-site constant up to the factor $(1-r)^{-1}$, which tends to $1$ when the dependence vanishes: this is what we call a dependent tensorization. The interaction only enters through the interdependence coefficients, which are Lipschitz constants measured with transport distances. For the Coulomb gas, we show that these coefficients only involve integrated quantities of the interaction: $\|W\|_{L^1}$ for $c^{\delta_0,1}_{ij}$ (Lemma~\ref{DIC}) and the $H^{-1}(\mu^{\y})$-norm of $\nabla W$ for $c^{d,2}_{ij}$ (Proposition~\ref{P-coef_W2}).

\subsection{Dobrushin interdependence coefficients and the uniform Poincar\'e inequality}\label{sec:tv}

Recall that for probability measures $\mu$ and $\nu$, the total variation distance is
\bbeq\label{tv}
\|\mu-\nu\|_{\rm tv}=\sup_{A}|\mu(A)-\nu(A)|=\sup_{\osc(f) \le  1} |\mu(f)-\nu(f)|
=\frac{1}{2} \sup_{|f|\le  1} |\mu(f)-\nu(f)|.
\nneq
By the Kantorovich-Rubinstein duality (or the maximal coupling), $\|\mu-\nu\|_{\rm tv}=W_{1,\delta_0}(\mu,\nu)$ for the discrete metric $\delta_0$.

\bdef For $i\neq j$, the Dobrushin interdependence coefficient of $x^i$ on $x^j$ in total variation is
\bbeq
c_{ij}^{\scriptscriptstyle \rm tv} := \sup_{x^{-i}=y^{-i} \text{ off }j} \|\mu_i(\cdot|x^{-i})- \mu_i(\cdot|y^{-i})\|_{\rm tv},\qquad c^{\scriptscriptstyle \rm tv}_{ii}:=0 .
\nneq
\ndef
By \eqref{tv}, $c_{ij}^{\scriptscriptstyle \rm tv}=c_{ij}^{\delta_0,1}$. In words, $c_{ij}^{\rm tv}$ measures the maximal influence, in total variation, that a displacement of particle $j$ has on the conditional distribution of particle $i$, all other particles being fixed. The matrix $C^{\rm tv}:=(c_{ij}^{\scriptscriptstyle \rm tv})$ is the Dobrushin interdependence matrix. Following \cite[Theorem~2, condition (2.5)]{Dobrushin68}, the \emph{Dobrushin uniqueness condition} holds if $\max_{i}\sum_{j} c_{ij}^{\scriptscriptstyle \rm tv}<1$. Here $c_{ij}^{\scriptscriptstyle \rm tv}$ is Dobrushin's coefficient \cite[(2.4)]{Dobrushin68}, with the same normalization \eqref{tv} \cite[(2.3) and (2.20)]{Dobrushin68}, see also \cite{Dobrushin70}. Since $C^{\rm tv}$ has nonnegative entries, this condition implies that its spectral radius satisfies $r_{\rm sp}(C^{\rm tv})<1$.

\blem\label{DIC} For all $i\neq j$,
\bbeq\label{DIC0}
c_{ij}^{\scriptscriptstyle \rm tv} \le  \frac{\Gamma}{N-1}\, e^{-\Gamma W_*} \|W\|_{L^1(\T^3, dx)}.
\nneq
\nlem
An estimate of this type was obtained by Becker and Menegaki \cite[Proposition~3.1]{Becker-Menegaki26} for Coulomb gases on $\rr^d$, by bounding $\|\mu_i(\cdot|x^{-i})- \mu_i(\cdot|y^{-i})\|_{\rm tv}$ by the sum of the distances of both measures to the conditional distribution without the $j$-th pole. We instead follow the interpolation method of \cite[Proposition~5.5]{Wu_AOP06}, which estimates the variation along a curve from $\mu_i(\cdot|x^{-i})$ to $\mu_i(\cdot|y^{-i})$. Our proof gives, for every $x^{-i}=y^{-i}$ off $j$,
\bbeq\label{DIC_fine}
\|\mu_i(\cdot|x^{-i})-\mu_i(\cdot|y^{-i})\|_{\rm tv}\le\frac{\Gamma}{2(N-1)}\,e^{-\Gamma W_*}\int_{\T^3}\big|W(x-x^j)-W(x-y^j)\big|\,dx .
\nneq
The right-hand side tends to $0$ as $y^j\to x^j$, and \eqref{DIC_fine} holds for every value of the charges, since only $Z_t\ge1$ is used below.

\bprf
Let $x^{-i}=y^{-i}$ off $j$. For $x\in\T^3$ and $t\in[0,1]$, let
$V_t(x):= \frac{\Gamma}{N-1} \big[\sum_{k\neq i, j} W(x-x^k) +  (1-t)W(x-x^j)+t W(x-y^j)\big]$
and
$\mu_t(dx) =\rho_t(x)\,dx = Z_t^{-1}e^{-V_t(x)}dx$, where, by Jensen's inequality and $\int V_t\,dx=0$, $Z_t:=\int e^{-V_t}\,dx\ge 1$. The curve $t\mapsto\mu_t$ goes from $\mu_0=\mu_i(\cdot|x^{-i})$ to $\mu_1=\mu_i(\cdot|y^{-i})$. Because $Z_t\ge 1$, and since the coefficients of $W$ in $V_t$ are nonnegative and sum to $\Gamma$, we have $V_t\geqslant \Gamma W_*$ and $\rho_t\le  e^{-\Gamma W_*}$.
Let $\hat W(x):=W(x-x^j)-W(x-y^j)$, which satisfies $\int \hat W\,dx=0$ and $-\partial_tV_t=\frac{\Gamma}{N-1}\hat W$. Let $f$ be bounded measurable with $a\le  f\le  b$ and $b-a\le  1$. Because $|\hat W|\rho_t\le  e^{-\Gamma W_*}|\hat W|$ is integrable, we may differentiate under the integral, and we get that:
\bbeq\label{DIC1}
\mu_1(f)-\mu_0(f)= \frac{\Gamma}{N-1}\int_{0}^{1} \cov_{\mu_t}(f, \hat W)\, dt.
\nneq
Since $0\le \rho_t\le  e^{-\Gamma W_*}$ and $a-\mu_t(f)\le 0\le  b-\mu_t(f)$, we have $e^{-\Gamma W_*}(a- \mu_t(f)) \le  \rho_t(f-\mu_t(f)) \le  e^{-\Gamma W_*}(b - \mu_t(f))$. With the midpoint $m_t:=e^{-\Gamma W_*}( \frac{a+b}{2} -\mu_t(f))$ of this interval, we obtain
$\left|\rho_t(x)\big(f(x)-\mu_t(f)\big) - m_t \right| \le  \frac12 e^{-\Gamma W_*}$.
Using $\int \hat W\,dx=0$, we deduce that:
\begin{align}
\nonumber
 |\cov_{\mu_t}(f, \hat W)| &=\left|\int_{\T^3} \Big[\rho_t(x)\big(f(x)-\mu_t(f)\big) - m_t \Big] \hat W(x)\,dx\right|\\
 \label{DIC4}
 &\le  \frac12 e^{-\Gamma W_*}  \int_{\T^3} |W(x-x^j)-W(x-y^j)|\,dx\le  e^{-\Gamma W_*}  \|W\|_{L^1(\T^3,dx)},
\end{align}
where we used the triangle inequality and the translation invariance of $dx$. Plugging \eqref{DIC4} into \eqref{DIC1} and taking the supremum over all such $f$ gives \eqref{DIC_fine} and \eqref{DIC0}.
\nprf

\bprf[Proof of Theorem \ref{T_W2I}\textit{(1)}]
Under \eqref{T_W2I1}, Lemma~\ref{DIC} gives
$\max_{i}\sum_{j} c_{ij}^{\scriptscriptstyle \rm tv} \le  r_{\rm D}<1$. Since $C^{\rm tv}$ has nonnegative entries, its spectral radius is bounded by its maximal row sum, i.e.\ $r_{\rm sp}(C^{\rm tv})\le r_{\rm D}<1$. Moreover, each conditional distribution $\mu_i(\cdot|x^{-i})$ is a one-particle measure \eqref{mu_single} with total charge $\bar\gamma=\Gamma<2\pi$. By Proposition~\ref{PI_single}, it satisfies the Poincar\'e inequality with the constant $C_{{\rm P},1}=\frac{1}{2\pi(2\pi-\Gamma)}$, uniformly in $i$ and $x^{-i}$. We apply \cite[Theorem~2.1, (2.4)]{Wu_AOP06}, which holds for the coefficients of any metric, in particular for the discrete one. It gives
$$
(1-r_{\rm D})\var_{\mu_N}(f)\le\int_{(\T^3)^N}\sum_{i=1}^N\var_{\mu_i(\cdot|x^{-i})}(f)\,d\mu_N(x).
$$
Bounding each conditional variance by Proposition~\ref{PI_single}, applied to $x^i\mapsto f(x)$ with $x^{-i}$ fixed, and using $\sum_{i}|\nabla_{i}f|^2=|\nabla f|^2$, we get that:
$
(1-r_{\rm D})\var_{\mu_N}(f)\le C_{{\rm P},1}\int|\nabla f|^2\,d\mu_N ,
$
which is \eqref{T_PI}. This is the argument of \cite[Theorem~2.2]{Wu_AOP06}, and it proves Theorem~\ref{T-abstract-W}(a) in our setting.
\nprf

\subsection{$W_2$-interdependence coefficients and the uniform $W_2I$ inequality}\label{sec:W2coef}

\subsubsection{The $W_2$-decomposition inequality}
For $i\neq j$, the $W_2$-interdependence coefficient of $x^i$ on $x^j$ is
\bbeq\label{def_cij2}
\cij := c_{ij}^{d,2}=\sup_{x^{-i}=y^{-i} \text{ off }j} \frac{W_{2}\big(\mu_i(\cdot|x^{-i}), \mu_i(\cdot|y^{-i})\big)}{d(x^j, y^j)},\qquad c^{(2)}_{ii}:=0,
\nneq
and $C^{(2)}:=(\cij)_{1\le i,j\le N}$ is the $W_2$-interdependence matrix, with operator norm $\|C^{(2)}\|_{\rm op}$ on $\ell^2(\{1,\ldots,N\})$.

\begin{remark}\label{rem:cij2}
$\,$

(i) Replacing $W_2$ by $W_1$ in \eqref{def_cij2} gives the interdependence coefficient of \cite[(2.3)]{Wu_AOP06} for the metric $d$, which is bounded by $\cij$ since $W_1\le W_2$.

(ii) The tensorization inequality \cite[Theorem~2.1, (2.4)]{Wu_AOP06} for the variance holds for the coefficients of any metric on $\T^3$, since the variance does not involve the metric. This is why we could use the discrete metric in Section~\ref{sec:tv}. On the contrary, the $W_2I$ inequality is formulated with the metric $d$, and the proof of Theorem~\ref{T-GW} below requires that a displacement of $x^j$ by $d(x^j,y^j)$ moves $\mu_i(\cdot|x^{-i})$ by at most $\cij\, d(x^j,y^j)$ in $W_2$. A total variation bound does not provide such a linear control, and the coefficients $c_{ij}^{\scriptscriptstyle \rm tv}$ are not sufficient to tensorize the $W_2I$ inequality.

(iii) The condition $\|C^{(2)}\|_{\rm op}<1$ plays the role of the Dobrushin uniqueness condition with respect to $W_2$. Because $\mu_N$ is exchangeable, all the off-diagonal coefficients are equal, $\cij=c^{(2)}$, and $C^{(2)}=c^{(2)}(J-\mathrm{Id})$, where $J$ is the matrix with all entries equal to $1$. Its eigenvalues are $(N-1)c^{(2)}$ and $-c^{(2)}$, and
$
\|C^{(2)}\|_{\rm op}=\max_{i}\sum_{j\neq i}\cij=(N-1)\,c^{(2)} .
$
\end{remark}

Let $W_2$ be the Wasserstein distance on $\mathcal P_1((\T^3)^N)$ associated with the metric $d_N$. The following result is contained in the unpublished preprint \cite{GW06}.

\bthm\label{T-GW} {\rm (\cite[Theorem~3.1]{GW06})} If $\|C^{(2)}\|_{\rm op}<1$, then for all $\nu\in \mathcal P_1((\T^3)^{N})$,
\bbeq \label{T-GW1}
W_{2}(\nu,\mu_N)^2 \le  \frac{1}{(1- \|C^{(2)}\|_{\rm op})^2} \int_{(\T^3)^N}
\sum_{i=1}^N W_{2}\big(\nu_i(\cdot|y^{-i}),\,\mu_i(\cdot|y^{-i})\big)^2\,\nu(dy),
\nneq
where $\nu_i(\cdot|y^{-i})$ denotes a regular version of the conditional distribution\footnote{I.e.\ a probability kernel $\nu_i(dx^i|x^{-i})$ from $(\T^3)^{N-1}$ to $\T^3$ such that $\nu(dx^i,dx^{-i})=\nu_i(dx^i|x^{-i})\,\nu^{-i}(dx^{-i})$, where $\nu^{-i}$ is the marginal law of $x^{-i}$ under $\nu$. Since $\T^3$ is Polish, such a kernel exists and is unique up to a $\nu^{-i}$-null set \cite[Theorems~5.3 and~5.4]{Kallenberg97}. In particular, the right-hand side of \eqref{T-GW1} does not depend on the choice of the version.} of $x^i$ under $\nu$ knowing $x^{-i}=y^{-i}$.
\nthm

For the reader's convenience, we give the proof of Theorem~\ref{T-GW}, taken from \cite{GW06}. It works on any product $E^N$ of a Polish metric space $(E,d)$, equipped with the $\ell^2$-product metric.

\bprf
Write $\mu=\mu_N$. Let $(X,Y)$ be a pair of random variables whose law is an optimal coupling of $(\mu,\nu)$ for $W_2$, i.e.
\begin{equation}\label{eq.optC}
\ee \Big [\sum_{i=1}^N  d(X^i, Y^i)^2 \Big ]= W_{2}(\mu,\nu)^2.
\end{equation}
Fix $i$, and let $\tilde X$, $\tilde Y$ be such that $\tilde X^j= X^j$ and $\tilde Y^j= Y^j$ for all $j\neq i$, and such that the conditional law of $(\tilde X^i, \tilde Y^i)$ given $(X^{-i}, Y^{-i})$ is an optimal coupling of $\big(\mu_i(\cdot|X^{-i}), \nu_i(\cdot|Y^{-i})\big)$ for $W_2$. By \cite[Corollary~5.22]{Villani09}, applied with the cost $d^2$ and the measurable map $(x^{-i},y^{-i})\mapsto(\mu_i(\cdot|x^{-i}),\nu_i(\cdot|y^{-i}))$, this optimal coupling can be chosen to depend measurably on $(X^{-i},Y^{-i})$. Since this conditional law has marginals $\mu_i(\cdot|X^{-i})$ and $\nu_i(\cdot|Y^{-i})$, which depend only on $X^{-i}$ and $Y^{-i}$ respectively, we have
\bbeq\label{T-GW-laws}
\tilde X\sim\mu,\qquad \tilde Y\sim\nu .
\nneq
Let us check it for $\tilde X$, assuming without loss of generality that $i=1$. Let $m$ be the law of $x^{-1}$ under $\mu$, let $F$ be bounded measurable, and set $\Phi(x^{-1}):=\int F(u,x^{-1})\,\mu_1(du\,|\,x^{-1})$. Then
$$
\ee\big[F(\tilde X^1,X^{-1})\big]
 \overset{(1)}{=}\ee\Big[\,\ee\big[F(\tilde X^1,X^{-1})\,\big|\,X^{-1},Y^{-1}\big]\Big]
\overset{(2)}{=}\ee\big[\Phi(X^{-1})\big]
\overset{(3)}{=}\int\Phi\,dm
\overset{(4)}{=}\int F\,d\mu .
$$
Here (1) is the tower property; (2) holds because, given $(X^{-1},Y^{-1})$, $\tilde X^1$ has law $\mu_1(\cdot|X^{-1})$; (3) holds because $X^{-1}$ has law $m$; and (4) is the disintegration $\mu(du,dx^{-1})=\mu_1(du|x^{-1})\,m(dx^{-1})$. The same argument gives $\tilde Y\sim\nu$. As a result, the law of $(\tilde X, \tilde Y)$ is again a coupling of $(\mu,\nu)$, and by optimality
$$
\ee \Big[\sum_{k=1}^N d(X^k, Y^k)^2\Big] \le  \ee \Big[\sum_{k=1}^N d(\tilde X^k,\tilde Y^k)^2\Big]
= \ee \Big[\sum_{k\ne i} d(X^k, Y^k)^2\Big] + \ee\big[d(\tilde X^i,\tilde Y^i)^2\big],
$$
which entails
\begin{equation}\label{eq.dxx}
\ee\big[d(X^i, Y^i)^2\big] \le   \ee\big[d(\tilde X^i,\tilde Y^i)^{2}\big] = \ee \Big[W_{2}\big(\mu_i(\cdot|X^{-i}), \nu_i(\cdot|Y^{-i})\big)^2\Big].
\end{equation}
On the other hand, by the triangle inequality,
$$
W_{2} \big(\mu_i(\cdot|X^{-i}), \nu_i(\cdot|Y^{-i})\big)\le W_{2} \big(\mu_i(\cdot|X^{-i}), \mu_i(\cdot|Y^{-i})\big)+ W_{2} \big(\mu_i(\cdot|Y^{-i}),\nu_i(\cdot|Y^{-i})\big).
$$
To bound the first term, we change the coordinates $j\neq i$ one at a time from $X^j$ to $Y^j$. Two consecutive configurations of this chain differ only in one coordinate $j\neq i$, and by the definition of $\cij$, the $W_2$-distance between the corresponding conditional laws is at most $\cij\,d(X^j,Y^j)$. By the triangle inequality along the chain, we get that:
$$
W_{2} \big(\mu_i(\cdot|X^{-i}), \nu_i(\cdot|Y^{-i})\big)\le  \sum_{j\ne i}
\cij\,  d(X^{j}, Y^{j})+ W_{2} \big(\mu_i(\cdot|Y^{-i}),\nu_i(\cdot|Y^{-i})\big).
$$
Squaring, taking expectations and using \eqref{eq.dxx}, we deduce that for each $i$,
\bbeq\label{T-GW3}
\ee\, [d(X^i,Y^i)^2] \le  \ee \Big[ \Big( \sum_{j\ne i} \cij\, d(X^{j}, Y^{j}) + W_{2} \big(\mu_i(\cdot|Y^{-i}), \nu_i(\cdot|Y^{-i})\big) \Big)^2 \Big].
\nneq
Let $U:=(d(X^i,Y^i))_{1\le i\le N}$ and $V:= \big(W_{2} (\mu_i(\cdot|Y^{-i}),\nu_i(\cdot|Y^{-i}))\big)_{1\le i\le N}$ be random column vectors, and for a random vector $Z$ in $\rr^N$, let $\|Z\|_2:=(\ee [\sum_{i} Z_i^2])^{1/2}$. Since $c^{(2)}_{ii}=0$, the right-hand side of \eqref{T-GW3} is $\ee[((C^{(2)}U+V)_i)^2]$. Summing \eqref{T-GW3} over $i$, and using Minkowski's inequality, we get
$
\|U\|_2 \le  \|C ^{(2)}U + V\|_2\le \|C^{(2)}\|_{\rm op}\|U\|_2+\|V\|_2 .
$
Since $\|U\|_2<\infty$ and $\|C^{(2)}\|_{\rm op}<1$, it follows that $\|U\|_2 \le  \frac {\|V\|_2}{1- \|C^{(2)}\|_{\rm op}}$. Finally, $\|U\|_2^2=W_2(\nu,\mu)^2$ by \eqref{eq.optC}, and since $Y$ has law $\nu$,
$\|V\|_2^2=\int\sum_{i} W_{2}\big(\nu_i(\cdot|y^{-i}),\mu_i(\cdot|y^{-i})\big)^2\,\nu(dy)$. This is \eqref{T-GW1}.
\nprf

\subsubsection{Inf-convolution and Hamilton-Jacobi equation}\label{sec.HJ}

We recall the relation between $W_2$, the inf-convolution and the Hamilton-Jacobi equation, see \cite[Sections~3.3.2 and~10.3.4]{Evans10} and \cite{Villani21}. For a Lipschitz function $f$ on $\T^3$, consider the inf-convolution
\bbeq\label{lem51d}
f_0=f, \qquad f_t(x)=Q_tf(x):=\inf_{y\in\T^3}\left\{f(y)+ \frac 1t\, d^2(x,y)\right\},\quad x\in\T^3,\ t>0.
\nneq

\blem\label{le.HJ}
Let $f$ be a Lipschitz function on $\T^3$ with Lipschitz constant $\|f\|_{\rm Lip}$. Then:
\begin{enumerate}[(a)]
  \item For every $x\in\T^3$, $t\mapsto Q_tf(x)$ is nonincreasing on $[0,\infty)$, and $\min f\le Q_tf\le f$.
  \item $(t,x)\mapsto Q_tf(x)$ is Lipschitz on $[0,\infty)\times\T^3$, and $\|Q_tf\|_{\rm Lip}\le\|f\|_{\rm Lip}$ for all $t\ge0$.
  \item $Q_tf$ satisfies the Hamilton-Jacobi equation
  \bbeq\label{lem51e}
  \partial_t f_t + \frac 14 |\nabla f_t|^2=0
  \nneq
  at every point of differentiability of $(t,x)\mapsto f_t(x)$, and thus for a.e.\ $(t,x)\in(0,\infty)\times\T^3$. In particular, for a.e.\ $t>0$, \eqref{lem51e} holds $dx$-a.e., and thus $\mu$-a.e.\ for every $\mu\ll dx$.
  \item {\rm (Kantorovich duality)} For all probability measures $\mu,\nu$ on $\T^3$,
  \bbeq\label{A32}
  W_{2}^2(\nu, \mu) = \sup\left\{ \int_{\T^3} Q_1 f\, d\nu - \int_{\T^3} f\, d\mu :\ f\in C^\infty(\T^3) \right\}.
  \nneq
  \item $(t,x)\mapsto Q_tf(x)$ is the unique bounded and uniformly continuous viscosity solution of \eqref{lem51e} on $(0,\infty)\times\T^3$ with initial datum $f$.
\end{enumerate}
\nlem
 
The proof of Lemma \ref{le.HJ} is made in the appendix (see Section \ref{sec.HJ2}). 
Note that $f_t$ is in general only Lipschitz for $t>0$, even when $f\in C^\infty(\T^3)$, which is why \eqref{lem51e} only holds almost everywhere. In the sequel we only use \emph{(a)}-\emph{(d)}.

\subsubsection{Estimate of the $W_2$-interdependence coefficients}

\bprop\label{P-coef_W2} Assume that $\Gamma<2\pi$. Then, for all $i\neq j$,
\bbeq\label{P-coef_W2b}
\cij  \le  \frac{\Gamma}{N-1}\, e^{-\Gamma W_*/2} \sqrt{ C_{{\rm So}} \big[1+(1+2\pi\Gamma )C_{{\rm P},1}\big]}\,\sup_{e\in \rr^3,\, |e|=1}\|e\cdot\nabla W \|_{L^{6/5}(\T^3,dx)}.
\nneq
\nprop

\bprf The strategy is taken from the unpublished preprint \cite{GW06}: we move the pole $x^j$ to $y^j$ along a geodesic, and we follow the Kantorovich potential along the Hamilton-Jacobi flow.

\noindent\textbf{Step 0. Setting.} Let $x^{-i},y^{-i}$ be such that $x^k=y^k$ for all $k\neq i,j$, and $x^j\neq y^j$. Let $\eta$ be a constant-speed geodesic in $\T^3$ from $\eta(0)=x^j$ to $\eta(1)=y^j$, with $|\dot{\eta}|=d(x^j,y^j)>0$, and let $e:=\dot{\eta}/d(x^j,y^j)$. For $t\in [0,1]$, let
$$
\aligned
\mu_t (dx)&=\rho_t(x)\,dx= \frac{1}{Z(t)}\exp\Big( -\frac{\Gamma }{N-1}\sum_{k\neq i,j} W(x-x^k) - \frac{\Gamma }{N-1}\, W(x-\eta(t)) \Big) dx, \\
\LL_t &= \Delta  -\Big(\frac{\Gamma }{N-1}\sum_{k\neq i,j} \nabla W(x-x^k) + \frac{\Gamma }{N-1}\, \nabla W(x-\eta(t))\Big)\cdot \nabla .
\endaligned
$$
Then, $\mu_0=\mu_i(\cdot|x^{-i})$ and $\mu_1=\mu_i(\cdot|y^{-i})$. Each $\mu_t$ is a one-particle measure \eqref{mu_single} with total charge $\Gamma$ and generator $\LL_t$, and $\rho_t\le e^{-\Gamma W_*}$ by \eqref{eq:Zbound}. Let $G_t:=\nabla W(\cdot-\eta(t))$, and let $\bu_t=(u_1,u_2,u_3)$ be the $\mu_t$-centered solution of the Poisson equations
$-\LL_t u_k = G_{t,k}-\mu_t(G_{t,k})$, $k=1,2,3$.
Finally, fix $f\in C^\infty(\T^3)$ and let $f_t=Q_tf$. By Lemma~\ref{le.HJ}(b), $\|f_t\|_{\rm Lip}\le\|f\|_{\rm Lip}$, which implies $f_t\in H^1(\mu_t)$ and $|f_{t}(x)-f_s(x)|\le\frac14\|f\|_{\rm Lip}^2|t-s|$ for all $x$ and $s,t\ge0$.

 \textbf{Step 1. Regularity of $t\mapsto\mu_t$.} For $x\neq\eta(\tau)$, the map $\tau\mapsto\rho_\tau(x)$ is $C^1$ and, since $-\partial_\tau\big[W(x-\eta(\tau))\big]=\dot\eta\cdot\nabla W(x-\eta(\tau))$,
$$
\partial_\tau\rho_\tau=\frac{\Gamma}{N-1}\,\rho_\tau\Big(\dot\eta\cdot G_\tau-\mu_\tau\big(\dot\eta\cdot G_\tau\big)\Big).
$$
Since $\rho_\tau\le e^{-\Gamma W_*}$ and $\int_{\T^3}|G_\tau|\,dx=\|\nabla W\|_{L^1(\T^3,dx)}<\infty$ (as $|\nabla W(x)|\sim|x|^{-2}$ near $0$), we have $\int_0^1\!\int|\partial_\tau\rho_\tau|\,dx\,d\tau<\infty$. By Fubini's theorem, for every bounded measurable $g$ and $s,t\in[0,1]$,
\bbeq\label{mut_int}
\mu_t(g)-\mu_s(g)=\frac{\Gamma}{N-1}\int_s^t\cov_{\mu_\tau}\big(g,\dot\eta\cdot G_\tau\big)\,d\tau .
\nneq
We give two consequences of \eqref{mut_int}.

(i) Because $\int_{\T^3}\dot\eta\cdot G_\tau\,dx=0$ and $\rho_\tau\le e^{-\Gamma W_*}$, the centering argument of the proof of Lemma~\ref{DIC} gives $|\cov_{\mu_\tau}(g,\dot\eta\cdot G_\tau)|\le\frac12e^{-\Gamma W_*}d(x^j,y^j)\|\nabla W\|_{L^1}$ when $\osc(g)\le1$. Thus, by \eqref{tv},
$
\|\mu_t-\mu_s\|_{\rm tv}\le\frac{\Gamma}{2(N-1)}\,e^{-\Gamma W_*}\,d(x^j,y^j)\,\|\nabla W\|_{L^1(\T^3,dx)}\,|t-s| .
$

(ii) If $g$ is continuous, then $\tau\mapsto\cov_{\mu_\tau}(g,\dot\eta\cdot G_\tau)$ is continuous. Indeed, by the change of variables $x=y+\eta(\tau)$,
$$
\mu_\tau\big(g\,\dot\eta\cdot G_\tau\big)=\int_{\T^3}g\big(y+\eta(\tau)\big)\,\rho_\tau\big(y+\eta(\tau)\big)\,\dot\eta\cdot\nabla W(y)\,dy,
$$
where $\rho_\tau\big(y+\eta(\tau)\big)= {Z(\tau)^{-1}}\,e^{-cW(y)}\prod_{k\neq i,j}e^{-cW(y+\eta(\tau)-x^k)}$ with $c=\frac{\Gamma}{N-1}$. The integrand is continuous in $\tau$ for every $y\neq0$ and bounded by $\|g\|_\infty e^{-\Gamma W_*}|\dot\eta|\,|\nabla W(y)|\in L^1(\T^3,dy)$, and dominated convergence applies. The same argument applies to $\mu_\tau(g)$ and $\mu_\tau(\dot\eta\cdot G_\tau)$. Consequently, by \eqref{mut_int}, $\frac{d}{dt}\mu_t(g)=\frac{\Gamma}{N-1}\cov_{\mu_t}(g,\dot\eta\cdot G_t)$ for every $t$.

 \textbf{Step 2. Derivative of $t\mapsto\mu_t(f_t)$.} By Step~0 and Step~1(i), the map $t\mapsto\mu_t(f_t)$ is Lipschitz (write $|\mu_t(f_t)-\mu_s(f_s)|\le |\mu_t(f_t-f_s)|+|\mu_t(f_s)-\mu_s(f_s)|$), and thus absolutely continuous on $[0,1]$. For $h\neq0$ small, we write
$$
\frac{\mu_{t+h}(f_{t+h})-\mu_t(f_t)}{h}
= \underbrace{\int_{\T^3}q_h\,\rho_{t}\,dx+\int_{\T^3}q_h\,(\rho_{t+h}-\rho_t)\,dx}_{=:A_h}
+ \underbrace{\frac{\mu_{t+h}(f_t)-\mu_t(f_t)}{h}}_{=:B_h},
$$
where $q_h:= ({f_{t+h}-f_t})/{h}$. By Step~0, $|q_h|\le\frac14\|f\|_{\rm Lip}^2$, and the second integral in $A_h$ is bounded by $\frac12\|f\|_{\rm Lip}^2\|\mu_{t+h}-\mu_t\|_{\rm tv}$, which tends to $0$ by Step~1(i). Since $(t,x)\mapsto f_t(x)$ is Lipschitz, it is differentiable at a.e.\ $(t,x)$ by Rademacher's theorem, and by Fubini's theorem, for a.e.\ $t$, $q_h\to\partial_tf_t$ $dx$-a.e. as $h\to 0$. By dominated convergence, the first integral in $A_h$ tends to $\mu_t(\partial_tf_t)$. Moreover, $B_h\to\frac{\Gamma}{N-1}\cov_{\mu_t}(\dot\eta\cdot G_t,f_t)$ by Step~1(ii) applied to $g=f_t$. Finally, by Lemma~\ref{le.HJ}(c), $\partial_tf_t=-\frac14|\nabla f_t|^2$ $\mu_t$-a.e.\ for a.e.\ $t$. Therefore, for a.e.\ $t\in[0,1]$,
\begin{equation}\label{eq.mut}
\frac d{dt} \mu_t (f_t)= \frac{\Gamma} {N-1}\cov_{\mu_t}\big( \dot{\eta}\cdot G_t, f_t\big) - \frac 14\, \mu_t(|\nabla f_t|^2).
\end{equation}

 \textbf{Step 3. Estimate of the derivative.} Since $\dot\eta=d(x^j,y^j)\,e$ and $\cov_{\mu_t}(e\cdot G_t,f_t)= \langle e\cdot G_t-\mu_t(e\cdot G_t),\, f_t \rangle_{\mu_t}$, we deduce from \eqref{eq.mut} that, for a.e.\ $t\in[0,1]$:
$$
\aligned \frac d{dt} \mu_t (f_t)
&= \frac{\Gamma\, d(x^j,y^j)} {N-1} \big\< -\LL_t (e\cdot\bu_t),\, f_t\big\>_{\mu_t}   - \frac 14\, \mu_t(|\nabla f_t|^2) \\
&=\frac{\Gamma\, d(x^j,y^j)} {N-1} \int_{\T^3} \nabla (e\cdot \bu_t) \cdot \nabla f_t\, d\mu_t- \frac 14\, \mu_t(|\nabla f_t|^2) \\
&\le  \left(\frac{\Gamma\, d(x^j,y^j)} {N-1}\right)^2 \int_{\T^3} |\nabla (e\cdot \bu_t)|^2\, d\mu_t
\le  \left(\frac{\Gamma\, d(x^j,y^j)} {N-1}\right)^2 K,
\endaligned$$
with $K:=e^{-\Gamma W_*} C_{{\rm So}}\big[1+(1+2\pi\Gamma)C_{{\rm P},1}\big]\sup_{|e|=1}\|e\cdot\nabla W\|^2_{L^{6/5}(\T^3,dx)}$. Here we used Proposition~\ref{prop:Wu}\emph{(c)} to integrate by parts against the Lipschitz function $f_t$, the inequality $ab-\frac14b^2\le  a^2$, and Lemma~\ref{lem26} applied to $\mu_t$, uniformly in $t$.

 \textbf{Step 4. Conclusion.} Integrating over $t\in[0,1]$, which is legitimate since $t\mapsto\mu_t(f_t)$ is absolutely continuous, we get
$$
\mu_i(Q_1f|y^{-i})- \mu_i(f|x^{-i})= \mu_1(f_1)-\mu_0(f_0) \le   \left(\frac{\Gamma\, d(x^j,y^j)} {N-1}\right)^2 K .
$$
Taking the supremum over all $f\in C^\infty(\T^3)$, we obtain by Lemma~\ref{le.HJ}(d) that
$
W_{2}^2 (\mu_0, \mu_1) \le  \big(\frac{\Gamma\, d(x^j,y^j)} {N-1}\big)^2 K .
$
Dividing by $d(x^j,y^j)^2$ and taking the supremum gives \eqref{P-coef_W2b}.
\nprf

\subsubsection{Proof of Theorem~\ref{T_W2I}\emph{(2)}}\label{sec:W2I}
Let $\nu\in \mathcal P_1((\T^3)^N)$. If $\nu=h^2\mu_N$ with $h\in H^1(\mu_N)$, then $\nu_i(\cdot|y)=h^2\mu_i(\cdot|y)/\mu_i(h^2|y)$ for $\nu$-a.e.\ $y$, and $I(\nu_i(\cdot|y)|\mu_i(\cdot|y))=\mu_i(|\nabla_ih|^2|y)/\mu_i(h^2|y)$. Since the law of $y^{-i}$ under $\nu$ has density $\mu_i(h^2|y)$ with respect to its law under $\mu_N$, we get the additivity of the Fisher-Donsker-Varadhan information \cite[Lemma~2.12]{GLWY_PTRF09}:
\bbeq\label{eq:additivity}
\int_{(\T^3)^N}\sum_{i=1}^N I\big(\nu_i(\cdot|y^{-i})\,|\,\mu_i(\cdot|y^{-i})\big)\,\nu(dy)=I(\nu|\mu_N).
\nneq

\bprf[Proof of Theorem \ref{T_W2I}\emph{(2)}]
By Proposition~\ref{P-coef_W2}, Remark~\ref{rem:cij2}\emph{(iii)}, Lemma~\ref{le.Sobolev} and the value of $C_{{\rm P},1}$, we have $\|C^{(2)}\|_{\rm op}=(N-1)c^{(2)}\le  r_{W_2}<1$. Let $\nu\in \mathcal P_1((\T^3)^N)$ with $I(\nu|\mu_N)<+\infty$. By Theorem~\ref{T-GW}, Proposition~\ref{W2I_1particle} applied to each $\mu_i(\cdot|y^{-i})$ (whose total charge is $\Gamma<2\pi$), and \eqref{eq:additivity}, we deduce that:
$$
\aligned
W_{2}(\nu,\mu_N)^2 &\le  \frac{1}{(1- r_{W_2})^2} \int_{(\T^3)^N} \sum_{i=1}^N W_{2}\big(\nu_i(\cdot|y^{-i}),\, \mu_i(\cdot|y^{-i})\big)^2\, \nu(dy)\\
&\le  \frac{(2 C_{W_2I,1})^2}{(1- r_{W_2})^2}  \int_{(\T^3)^N} \sum_{i=1}^N I\big(\nu_i(\cdot|y^{-i})\,|\,\mu_i(\cdot|y^{-i})\big)\,\nu(dy)
= \left(\frac{2 C_{W_2I,1}}{1- r_{W_2}}\right)^2  I(\nu|\mu_N).
\endaligned
$$
This is \eqref{T_W2I2}, and it proves Theorem~\ref{T-abstract-W}(b) in our setting.
\nprf

\subsubsection{Consequences}
\label{sec.cons}
By \cite[Theorem~2.4]{GLWW09}, the $W_2I$ inequality implies Talagrand's $W_2H$ inequality with the same constant.

\bcor\label{cor1}
Assume that $\Gamma<2\pi$ and $r_{W_2}<1$, and let $C_{W_2I}$ be given by \eqref{T_W2I5}. Then
\bbeq\label{T2}
 W_2^2(\nu, \mu_N) \le   2 C_{W_2I}\, H(\nu|\mu_N), \qquad \forall\, \nu\in \mathcal P_1 ((\T^3)^N).
\nneq
\ncor

As a consequence, $\mu_N$ satisfies a Gaussian concentration inequality which is dimension-free, i.e.\ uniform in $N$. Indeed, because $W_1\le W_2$, \eqref{T2} implies $W_1^2(\nu,\mu_N)\le 2C_{W_2I}\,H(\nu|\mu_N)$, which is equivalent, by Bobkov and G\"otze \cite{BobkovGotze99}, to $\int e^{\lambda(F-\mu_N(F))}d\mu_N\le e^{C_{W_2I}\lambda^2/2}$ for all $\lambda\in\rr$ and all functions $F$ which are $1$-Lipschitz with respect to $d_N$. In particular,
$
\mu_N\big(F-\mu_N(F)\geqslant r\big)\le e^{-r^2/(2C_{W_2I})}$ for $r>0$.
Moreover, the $W_2I$ inequality \eqref{T_W2I2} itself yields Gaussian deviation inequalities for the empirical time averages $\frac1t\int_0^tF(X_s)\,ds$ of the Langevin dynamics \eqref{Langevin_D}, uniformly in $N$, see \cite{GLWY_PTRF09}.

The second application is a lower bound on the ground state energy of the Schr\"odinger operator $-\LL + v$ on $L^2(\mu_N)$. For a function $v$ on $(\T^3)^N$, let
$Q_tv (x)=\inf_{y} \left( v(y) + \frac{1}{t}\, d_N^2(x,y)\right)$
be its inf-convolution, which solves the Hamilton-Jacobi equation $\partial_t Q_t v + \frac{1}{4}|\nabla Q_tv|^2=0$ with initial datum $v$ (see Lemma~\ref{le.HJ}, whose proof holds on $(\T^3)^N$). By \cite[Corollary~2.5, (2.a)$\Leftrightarrow$(2.c)]{GLWY_PTRF09}, we obtain the following result.

\bcor\label{cor2}
Assume that $\Gamma<2\pi$ and $r_{W_2}<1$, and let $c=C_{W_2I}$ be given by \eqref{T_W2I5}. Then, for every continuous function $v$ on $(\T^3)^N$,
\bbeq\label{cor2b}
\aligned
\lambda_0\left(-\LL + \frac{1}{4c^2}\, v \right)&=\inf\left\{ \frac{1}{4c^2}\int v\, h^2\, d\mu_N + \EE_{\mu_N}[h]:\ h\in H^1(\mu_N),\ \mu_N(h^2)=1\right\}\\
&\geqslant \frac{1}{4c^2} \int Q_1v\,d\mu_N.
\endaligned
\nneq
\ncor
In fact, the validity of \eqref{cor2b} for all such $v$ is equivalent to the $W_2I$ inequality \eqref{T_W2I2}, which is one of the main motivations for studying $W_2I$ inequalities.

\subsection{Extensions: the two-dimensional Coulomb gas and periodic Riesz gases}\label{sec:ext2}

The proof of Theorem~\ref{T_W2I} uses only a few properties of the kernel $W$. We identify them and deduce two extensions.

\subsubsection{The two-dimensional Coulomb gas}
Let $W$ be the periodic Green function on $\T^2$,
\bbeq\label{eq:W2D}
-\Delta W=2\pi(\delta_0-1),\qquad \int_{\T^2}W\,dx=0,
\nneq
i.e.\ $\hat W(k)=\frac{1}{2\pi|k|^2}$ for $k\neq0$. Then, $W$ is smooth on $\T^2\setminus\{0\}$, and $R(x):=W(x)+\log d(x,0)$ is bounded on $\T^2\setminus\{0\}$ and smooth near $0$ \cite{Serfaty26}. In particular $W_*:=\min_{\T^2}W>-\infty$, and $|\nabla W(x)|\le C\,d(x,0)^{-1}$, which implies $\nabla W\in L^p(\T^2,dx)$ for every $p<2$. Let $\mu_N$ be the log-gas \eqref{Coulomb-measure} on $(\T^2)^N$ with this $W$.


\bthm\label{T_2D} Assume that $\Gamma<4\pi$.
\begin{enumerate}
  \item[(1)] If $r^{\scriptscriptstyle (2D)}_{\rm D}:=\Gamma e^{-\Gamma W_*}\|W\|_{L^1(\T^2,dx)}<1$, then, for all $f\in H^1(\mu_N)$,
  $$
  \var_{\mu_N}(f)\le\frac{1}{(1-r^{\scriptscriptstyle (2D)}_{\rm D})\cdot\pi(4\pi-\Gamma)}\int_{(\T^2)^N}|\nabla f|^2\,d\mu_N .
  $$
  \item[(2)] If
  $$
  r^{\scriptscriptstyle (2D)}_{W_2}:=\Gamma e^{-\Gamma W_*/2}\sqrt{\Big(1+\frac{\sqrt2}{4}\Big)\Big[1+\frac{1+\pi\Gamma}{\pi(4\pi-\Gamma)}\Big]}\cdot\sup_{e\in\rr^2,\,|e|=1}\|e\cdot\nabla W\|_{L^{4/3}(\T^2,dx)}<1,
  $$
  then, for all $\nu\in\mathcal P_1((\T^2)^N)$, $W_2^2(\nu,\mu_N)\le[2C^{\scriptscriptstyle (2D)}_{W_2I}]^2\,I(\nu|\mu_N)$, with
  $$
  C^{\scriptscriptstyle (2D)}_{W_2I}=\frac{1}{1-r^{\scriptscriptstyle (2D)}_{W_2}}\sqrt{\Big(\frac14+\frac{24}{\pi(4\pi-\Gamma)}\Big)\frac1{\pi(4\pi-\Gamma)}} .
  $$
\end{enumerate}
\nthm

As in dimension $3$, the condition $\Gamma<4\pi$ is implied by the two other conditions: $\|W\|_{L^1(\T^2,dx)}\geqslant\hat W(1,0)=\frac1{2\pi}$ and $\|\partial_1W\|_{L^{4/3}(\T^2,dx)}\geqslant|\int\partial_1W(x)\sin(2\pi x_1)\,dx|=1$. With the numerical values of Section~\ref{sec.N}, the conditions of Theorem~\ref{T_2D}\emph{(1)} and~\emph{(2)} hold  for $\Gamma\le1.93$ and $\Gamma\le0.37$ respectively.

\bprf
We follow the proof of Theorem~\ref{T_W2I} and only detail the steps where the dimension enters. 

Let $V$ be the one-particle potential or the $N$-particle potential on $\T^n$, $n\in\{2,2N\}$, as in \eqref{eq.V1}-\eqref{eq.V2}. The singular set $\SS$ is a finite set of points if $n=2$ (the singular set of the one-particle potential $V^{\y}$ is the set of poles $\{y^2,\dots,y^N\}$). If $n=2N$, it is the finite union $\SS_N=\bigcup_{i<j}\{x^i=x^j\}$, where each $\{x^i=x^j\}$ is a subtorus of $\T^{2N}$ of dimension $2N-2$, i.e.\ of codimension $2$. In both cases $\mathcal D=\T^n\setminus\SS$ is open and connected, since removing finitely many closed submanifolds of codimension $2$ does not disconnect $\T^n$.
 Properties (\textbf{P2}) and (\textbf{P3}) hold (see Section \ref{subsec:FI}), with $\Delta W=2\pi$ on $\T^2\setminus\{0\}$.

 In Section~\ref{subsec:FI}, property (\textbf{P1}) is only used in Proposition~\ref{prop:Wu}, to apply \cite{Wu_PTRF99}. The assumptions (\textbf{H1}) and (\textbf{H2}) there only involve $\mathcal D$: $\varphi$ is continuous and strictly positive on $\mathcal D$ with $\nabla\varphi\in L^2_{\rm loc}(\mathcal D)$, and $|\nabla\varphi|^2/\varphi^2$ belongs to the local Kato class of $\mathcal D$. Both follow from (\textbf{P2}), since $\varphi$ is smooth on $\mathcal D$ and $|\nabla\varphi|^2/\varphi^2=\frac14|\nabla V|^2$ is locally bounded on $\mathcal D$. Only the finite energy condition \cite[(0.2)]{Wu_PTRF99} involves the singular set: the extension of $\varphi$ by $0$ on $\SS$ must be continuous on $\T^n$, with distributional gradient in $L^2(\T^n,dx)$. The smoothness of $e^{-V}$ in (\textbf{P1}) is thus not needed. In dimension $2$, $e^{-V}$ is not smooth, since near a pole it behaves like $d(x,y)^{a}$ for some $a>0$, and we replace (\textbf{P1}) by the weaker property:
\medskip

(\textbf{P1}') The function $\varphi:=Z^{-1/2}e^{-V/2}$  is continuous on $\T^n$, strictly positive on $\mathcal D$, and belongs to $H^1(\T^n,dx)$.
\medskip

Continuity and positivity are clear. We now prove   $\varphi\in H^1(\T^n,dx)$. 
\smallskip

\emph{The gradient on $\mathcal D$.} 
\smallskip

Recall that $R(z)=W(z)+\log d(z,0)$ is bounded, which gives $W(z)\geqslant-\log d(z,0)-\|R\|_\infty$ and $|\nabla W(z)|\le C\,d(z,0)^{-1}$ for $z\neq0$. Consider the $N$-particle potential \eqref{eq.V2},
$$V(x)=\sum_{i\neq j}\gamma_{ij}W(x^i-x^j).$$
Since $W$ is even, the pair $\{a,b\}$ has total coefficient $\gamma_{ab}+\gamma_{ba}$ in $V$ (for the log-gas, $\gamma_{ab}+\gamma_{ba}=\frac{\Gamma}{N-1}$).   Fix a pair $a\neq b$ and set $\gamma:=\gamma_{ab}+\gamma_{ba}$. The other terms of $V$ are bounded from below by multiples of $W_*$, and $Z\geqslant1$. Therefore, for all $x\in\T^{2N}$,
$$
\varphi(x)^2=\frac{e^{-V(x)}}{Z}\le C\,e^{-\gamma W(x^a-x^b)}\le C\,d(x^a,x^b)^{\gamma},
$$
where $C$ does not depend on $x$ (nor on the pair, by taking the maximum over the finitely many pairs). Moreover, $|\nabla V|^2\le C\sum_{a\neq b}d(x^a,x^b)^{-2}$ on $\mathcal D$. Since $\nabla\varphi=-\frac12\varphi\nabla V$ on $\mathcal D$, bounding $\varphi^2$ in each term by the factor of the corresponding pair, we get that:
$$
|\nabla\varphi|^2\le C\sum_{a\neq b}d(x^a,x^b)^{-2+\gamma_{ab}+\gamma_{ba}}\qquad\text{on }\mathcal D .
$$
Each term is integrable on $\T^{2N}$ (by the translation invariance of $dx$ and the use of polar coordinates). 
Thus the gradient of $\varphi$ on $\mathcal D$ belongs to $L^2(\T^n,dx)$. The one-particle case is identical, with the distinct poles $y$ (with total charges $a_y>0$) instead of the pairs.
\smallskip

\emph{The distributional gradient on $\T^n$.} 

\smallskip

The previous bound only concerns the gradient of $\varphi$ on $\mathcal D$. We now check that $\SS$ does not carry any singular part of the distributional gradient of $\varphi$ on $\T^n$ (as the derivative $\delta_0$ of a function jumping at $0$).
The argument   relies on two facts. First, $\varphi$ vanishes fast enough at $\SS$ for its gradient on $\mathcal D$ to be square integrable. Second, $\SS$ is small: since it has codimension $2$, there exist cutoff functions $\chi_\vep$ vanishing near $\SS$, equal to $1$ away from it, with $\int_{\T^n}|\nabla\chi_\vep|^2dx\to0$. This second fact is purely geometric and does not depend on $\varphi$.

 For $0<\vep<\frac12$, let $\psi_\vep(r)=0$ for $r\le\vep^2$, $\psi_\vep(r)=\log(r/\vep^2)/\log(1/\vep)$ for $\vep^2\le r\le\vep$, and $\psi_\vep(r)=1$ for $r\geqslant\vep$. the function $\psi_\vep$ is a cutoff function in the distance variable $r$ vanishing near $0^+$. It is   Lipschitz ($\|\psi_\varepsilon'\|_\infty
=1/({\varepsilon^2\log(1/\varepsilon)})<\infty$), but not $C^1$. Then, it holds that:
$$\int_{\rr^2}|\nabla[\psi_\vep(|z|)]|^2dz=2\pi\int_{\vep^2}^{\vep}\frac{dr}{r\log^2(1/\vep)}=\frac{2\pi}{\log(1/\vep)}\to0 .$$
Set $\chi_\vep(x):=\prod_{a<b}\psi_\vep(d(x^a,x^b))$ (resp.\ $\prod_y\psi_\vep(d(x,y))$ over the distinct poles). The function $\psi_\vep(d(x^a,x^b))$ is constant outside the region $\vep^2<d(x^a,x^b)<\vep$. Since $\vep<\frac12$, in this region the torus distance agrees locally with the Euclidean distance between suitable lifts of $x^a$ and $x^b$, it is smooth, and $|\nabla_{x^a}d(x^a,x^b)|=|\nabla_{x^b}d(x^a,x^b)|=1$. Consequently, the gradient of $\psi_\vep(d(x^a,x^b))$ with respect to $x$ has norm $\sqrt2\,|\psi_\vep'(d(x^a,x^b))|$ a.e.
 By translation invariance, its squared $L^2(\T^{2N},dx)$-norm equals $2\int_{\rr^2}|\nabla[\psi_\vep(|z|)]|^2dz$. Since $0\le\psi_\vep\le1$ and the number of pairs is finite, we deduce that when $\vep \to 0^+$:
 $$\int_{\T^n}|\nabla\chi_\vep|^2dx\to0.$$ The support of $\chi_\vep$ is contained in the compact set $\{d(x^a,x^b)\geqslant\vep^2\ \text{for all } a\neq b\}\subset\mathcal D$, on which $\varphi$ is smooth. As a consequence, $\chi_\vep\varphi$ is Lipschitz on $\T^n$, belongs to $H^1(\T^n,dx)$, and $\nabla(\chi_\vep\varphi)=\chi_\vep\nabla\varphi+\varphi\nabla\chi_\vep$. Moreover, $0\le\chi_\vep\le1$ and $\chi_\vep\to1$ pointwise on $\mathcal D$ when $\vep \to 0^+$, which has full measure. Recall that $\varphi$ is bounded on $\T^n$, being continuous on a compact set, and that $\nabla\varphi\in L^2(\T^n,dx)$ by the previous step, where $\nabla\varphi$ denotes the gradient computed on $\mathcal D$. Since
$$
|\chi_\vep\varphi-\varphi|\le|\varphi|\le\|\varphi\|_\infty\quad\text{and}\quad|\chi_\vep\nabla\varphi-\nabla\varphi|\le|\nabla\varphi| ,
$$
dominated convergence implies  $\chi_\vep\varphi\to\varphi$ and $\chi_\vep\nabla\varphi\to\nabla\varphi$ in $L^2(\T^n,dx)$, when $\vep \to 0^+$. In addition, $\|\varphi\nabla\chi_\vep\|_{2}\le\|\varphi\|_\infty\|\nabla\chi_\vep\|_2\to0$. Consequently, $\chi_\vep\varphi\to\varphi$ and $\nabla(\chi_\vep\varphi)=\chi_\vep\nabla\varphi+\varphi\nabla\chi_\vep\to\nabla\varphi$ in $L^2(\T^n,dx)$. Passing to the limit in $\int_{\T^n}\chi_\vep\varphi\,\partial_k\psi\,dx=-\int_{\T^n}\partial_k(\chi_\vep\varphi)\,\psi\,dx$, for $\psi\in C^\infty(\T^n)$ and $k\in\{1,\dots,n\}$, we get that:
$$
\int_{\T^n}\varphi\,\partial_k\psi\,dx=-\int_{\T^n}\partial_k\varphi\,\psi\,dx .
$$
Therefore, $\varphi\in H^1(\T^n,dx)$, and its distributional gradient on $\T^n$ is the gradient computed on $\mathcal D$. In particular, $\SS$ carries no singular part of this gradient.

\emph{Conclusion.} Assumptions (\textbf{H1}), (\textbf{H2}) and (0.2) of \cite{Wu_PTRF99} hold in dimension $2$. Since the proofs of Propositions~\ref{prop:wellposed} and~\ref{prop:Wu} only use (\textbf{P1}'), (\textbf{P2}) and (\textbf{P3}), these results hold in dimension $2$, and so does Proposition~\ref{prop:Wu}\emph{(c)}.

\emph{The other steps.} The rest of the proof follows the lines of the proof of Theorem~\ref{T_W2I}. We only indicate the changes in the constants.
\begin{itemize}
\item \emph{One particle.} The ground state computation gives $\int_{\T^2}|\nabla(f\varphi)|^2\,dx\le\EE_{\mu^{\y}}[f]+\pi\bar\gamma\int_{\T^2}f^2\,d\mu^{\y}$, and the Poincar\'e inequality of Proposition~\ref{PI_single} holds with $C^{\scriptscriptstyle (2D)}_{{\rm P},1}=\frac{1}{\pi(4\pi-\bar\gamma)}$ for $\bar\gamma<4\pi$. Lemma~\ref{le.Liu} applies with ${\rm Diam}(\T^2)^2=\frac12$, which gives the one-particle $W_2I$ inequality with constant $C^{\scriptscriptstyle (2D)}_{W_2I,1}=\sqrt{(\frac14+24C^{\scriptscriptstyle (2D)}_{{\rm P},1})C^{\scriptscriptstyle (2D)}_{{\rm P},1}}$.
\item \emph{Poincar\'e inequality.} The proof of Lemma~\ref{DIC} only uses $Z_t\geqslant1$, $V_t\geqslant\Gamma W_*$ and $W\in L^1$, and the proof of Theorem~\ref{T_W2I}\emph{(1)} gives~(1).
\item \emph{Sobolev inequality on $\T^2$.} For every $w\in H^1(\T^2,dx)$,
\bbeq\label{Sobolev2D}
\|w\|_{4}^2\le C^{\scriptscriptstyle (2D)}_{{\rm So}}\big(\|\nabla w\|^2_{2}+\|w\|^2_{2}\big),\qquad C^{\scriptscriptstyle (2D)}_{{\rm So}}\le1+\frac{\sqrt2}{4}.
\nneq
Indeed, as in the proof of Lemma~\ref{le.Sobolev}, for $f\in C^1(\T^2)$, $\|f\|_2\le\|f\|_1+\frac14(\|\partial_1f\|_1+\|\partial_2f\|_1)$. We apply this to $f=w^2$. Since $\|\partial_if\|_1\le2\|w\|_2\|\partial_iw\|_2$ and $\|\partial_1w\|_2+\|\partial_2w\|_2\le\sqrt2\,\|\nabla w\|_2$, we get that $\|w\|_4^2\le\|w\|_2^2+\frac{\sqrt2}{2}\|w\|_2\|\nabla w\|_2\le(1+\frac{\sqrt2}4)\|w\|_2^2+\frac{\sqrt2}4\|\nabla w\|_2^2$.
\item \emph{$H^{-1}$ estimate.} Near the pole $y^j$, the density of $\mu^{\y}$ is bounded by $C\,d(x,y^j)^{a}$ for some $a>0$, and $|\nabla W(x-y^j)|^2\le C\,d(x,y^j)^{-2}$. Since $\int_{|z|\le1}|z|^{a-2}\,dz<\infty$ in dimension $2$, $\nabla W(\cdot-y^j)\in L^2(\mu^{\y})$. The proof of Lemma~\ref{lem26}, with H\"older's inequality for the exponents $(4,\frac43)$ instead of $(6,\frac65)$ and \eqref{Sobolev2D} instead of \eqref{SobolevI}, gives
\bbeq\label{lem26_2D}
\int_{\T^2}|\nabla(e\cdot\bu^{\scriptscriptstyle (j)})|^2\,d\mu^{\y}\le e^{-\bar\gamma W_*}\,C^{\scriptscriptstyle (2D)}_{{\rm So}}\big[1+(1+\pi\bar\gamma)C^{\scriptscriptstyle (2D)}_{{\rm P},1}\big]\,\|e\cdot\nabla W\|^2_{L^{4/3}(\T^2,dx)} ,
\nneq
where $\|e\cdot\nabla W\|_{L^{4/3}(\T^2,dx)}<\infty$ since $|\nabla W|^{4/3}\le C\,d(x,0)^{-4/3}$ is integrable on $\T^2$.
\item \emph{$W_2I$ inequality.} Lemma~\ref{le.HJ} holds on $\T^2$ with the same proof. The proof of Proposition~\ref{P-coef_W2}, with \eqref{lem26_2D} instead of Lemma~\ref{lem26}, gives $\|C^{(2)}\|_{\rm op}\le r^{\scriptscriptstyle (2D)}_{W_2}$, and the end of the proof of Theorem~\ref{T_W2I}\emph{(2)} gives~(2).
\end{itemize}
\nprf

\subsubsection{Periodic Riesz gases}
For $0<s<3$, let $W_s$ be the $\zz^3$-periodic function with 
$$\int_{\T^3}W_s\,dx=0$$ and Fourier coefficients $\hat W_s(k)=c_s|k|^{s-3}$, $k\neq0$, where $c_s>0$ is such that $|x|^{-s}$ has Fourier transform $c_s|\xi|^{s-3}$ on $\rr^3$. Then, we have that:
 $W_s(x)=|x|^{-s}+R_s(x)$ near $0$, where $R_s$ is smooth, and $W_s$ is smooth on $\T^3\setminus\{0\}$. Let $\mu_N^{\scriptscriptstyle (s)}$ be the Gibbs measure \eqref{Coulomb-measure} with $W$ replaced by $W_s$.

The only place where the specific form of $W$ enters the one-particle Poincar\'e inequality is the ground state identity of Step~1 in the proof of Proposition~\ref{PI_single}, in which only an upper bound on $\Delta V$ away from the poles is needed. In dimension $3$, $\Delta|x|^{-s}=s(s-1)|x|^{-s-2}$ on $\rr^3\setminus\{0\}$, which is nonpositive for $0<s\le1$. As a result, $\Delta W_s$ is bounded from above on $\T^3\setminus\{0\}$, and we set $K_s:=\sup_{\T^3\setminus\{0\}}\Delta W_s<\infty$. Note that $K_s>0$, since $\int_{\T^3}\Delta W_s\,dx=0$, and that $K_1=4\pi$. Then, $\Delta V^{\y}\le\bar\gamma K_s$ off the poles, and the proof of Proposition~\ref{PI_single} gives the one-particle Poincar\'e inequality with constant $(4\pi^2-\frac12\bar\gamma K_s)^{-1}$, provided $\bar\gamma<8\pi^2/K_s$.

The other ingredients only use properties of $W$ that $W_s$ shares for $0<s\le1$. In Section~\ref{subsec:FI}, $W_s\to+\infty$ at $0$, and $e^{-V}$ vanishes with all its derivatives at the singular set, as in (\textbf{P1})-(\textbf{P2}); moreover $\LL V=\Delta V-|\nabla V|^2$ is bounded from above, which is the Lyapunov argument of Proposition~\ref{prop:wellposed}. In Lemma~\ref{DIC}, $\int W_s\,dx=0$, $W_{s,*}:=\min W_s>-\infty$, and $W_s\in L^1$. In Lemma~\ref{lem26}, $|\nabla W_s(x)|\sim s|x|^{-s-1}$ near $0$, and $\nabla W_s\in L^{6/5}(\T^3,dx)$ if and only if $s<\frac32$; furthermore, $|\nabla W_s|^2e^{-\gamma W_s}$ is bounded near $0$, which gives $\nabla W_s(\cdot-y^j)\in L^2(\mu^{\y})$. In Proposition~\ref{P-coef_W2}, $\nabla W_s\in L^1$. The tensorization arguments and Lemma~\ref{le.Liu} do not depend on the interaction. Consequently, we have the following theorem.

\bthm\label{T_Riesz}
Let $0<s\le1$, assume that $\Gamma<8\pi^2/K_s$, and set $C^{(s)}_{{\rm P},1}:=\big(4\pi^2-\frac12\Gamma K_s\big)^{-1}$.
\begin{enumerate}
  \item[(1)] If $r_{\rm D}^{(s)}:=\Gamma e^{-\Gamma W_{s,*}}\|W_s\|_{L^1(\T^3,dx)}<1$, then
  $
  \var_{\mu_N^{\scriptscriptstyle (s)}}(f)\le\frac{C^{(s)}_{{\rm P},1}}{1-r_{\rm D}^{(s)}}\int|\nabla f|^2\,d\mu_N^{\scriptscriptstyle (s)}$ for all $f\in H^1(\mu_N^{\scriptscriptstyle (s)})$.
  \item[(2)] If $r_{W_2}^{(s)}:=\Gamma e^{-\Gamma W_{s,*}/2}\sqrt{\frac73\big[1+(1+\frac12\Gamma K_s)C^{(s)}_{{\rm P},1}\big]}\,\sup_{|e|=1}\|e\cdot\nabla W_s\|_{L^{6/5}(\T^3,dx)}<1$, then
  $$
  W_2^2(\nu,\mu_N^{\scriptscriptstyle (s)})\le\Big(\frac{2C^{(s)}_{W_2I,1}}{1-r_{W_2}^{(s)}}\Big)^2I(\nu|\mu_N^{\scriptscriptstyle (s)}),\qquad\forall\,\nu\in\mathcal P_1((\T^3)^N),
  $$
  where $C^{(s)}_{W_2I,1}:=\sqrt{\big(\frac38+24C^{(s)}_{{\rm P},1}\big)C^{(s)}_{{\rm P},1}}$.
\end{enumerate}
For $s=1$, $K_1=4\pi$ and this is Theorem~\ref{T_W2I}.
\nthm


\subsubsection{Which properties of the interaction are used}\label{sec:discussionW}
Let us summarize the assumptions on the interaction used by the dependent tensorization. The one-site inequalities require $\int W\,dx=0$, $W$ bounded from below, and $\Delta W$ bounded from above away from $0$ (ground state transformation). The Dobrushin coefficients require $W\in L^1$. The $W_2$-interdependence coefficients require a bound on the $H^{-1}(\mu^{\y})$-norm of $\nabla W(\cdot-y^j)$, uniform in the poles. Our proof obtains it from the Sobolev embedding $H^1(\T^d)\subset L^{2d/(d-2)}(\T^d)$ ($d\geqslant3$), and thus requires $\nabla W\in L^{2d/(d+2)}(\T^d,dx)$. For a Riesz-type singularity $|x|^{-s}$, this holds if and only if $s<\frac d2$. For the Coulomb kernel $s=d-2$, this holds if and only if $d\le3$ (and the logarithmic case $d=2$ is covered by \eqref{lem26_2D}). Note that this condition only involves integrated quantities of $\nabla W$. This is the reason why the interdependence coefficients remain finite, although the force $\nabla W$ is singular.


\section{Tensorization with a defect: the uniform log-Sobolev inequality}\label{sec:LSI}

 In this section we prove Theorem~\ref{T_LSI_intro}. We first explain why the mechanism of Section~\ref{sec:W2I_all} does not extend to the entropy (Section~\ref{sec:obstruction}). We then prove a uniform one-site log-Sobolev inequality (Section~\ref{sec:LSI1}), establish an abstract propagation principle with a defect (Section~\ref{sec:defect}), and combine it with the partition bound of Wang and Zhao \cite{WangZhao26} (Sections~\ref{sec:partition} and~\ref{sec:LSI_proof}). Section~\ref{sec:LSI_ext} discusses extensions and the role of each assumption.

We set
\bbeq\label{eq:energy}
\mathbf E_N(x):=\frac{\Gamma}{2(N-1)}\sum_{{i,j=1, i\neq j}}^{N}W(x^i-x^j),\quad Z_N(t):=\int_{(\T^3)^N}e^{-t\mathbf E_N}\,dx^1\cdots dx^N,\quad t>0,
\nneq
so that $\mu_N=Z_N(1)^{-1}e^{-\mathbf E_N}dx^1\cdots dx^N$. We write $\mathbf m _N:=dx^1\cdots dx^N$ for the Haar probability measure on $(\T^3)^N$. By Jensen's inequality and $\int_{\T^3} W\,dx=0$, $Z_N(t)\geqslant1$ for all $t>0$.

\subsection{Why interdependence coefficients do not tensorize the entropy}\label{sec:obstruction}

The log-Sobolev inequality also satisfies the independent tensorization for product measures, and it is natural to look for a dependent tensorization in the spirit of Theorem~\ref{T-abstract-W}. The classical criteria are those of Otto and Reznikoff \cite{OttoReznikoff07}, Marton \cite{Marton13} and Zegarlinski \cite{Zegarlinski92}. In this subsection, we show that none of them applies to $\mu_N$, however small $\Gamma$ is, and we explain why the transport argument of Section~\ref{sec:W2I_all} does not meet the same obstruction.

\emph{The Hessian of the Hamiltonian.} The Hamiltonian of $\mu_N$ is $\mathbf E_N$, and for $i\neq j$,
$$
\nabla^2_{x^ix^j}\mathbf E_N(x)=-\frac{\Gamma}{N-1}\nabla^2W(x^i-x^j),\qquad
\nabla^2_{x^ix^i}\mathbf E_N(x)=\frac{\Gamma}{N-1}\sum_{k\neq i}\nabla^2W(x^i-x^k).
$$
Since $W(x)=|x|^{-1}+R(x)$ with $R$ smooth near $0$, and $\nabla^2|x|^{-1}=|x|^{-5}(3x\otimes x-|x|^2{\rm Id})$ has eigenvalues $2|x|^{-3}$ and $-|x|^{-3}$ (twice), the Hessian of $W$ is unbounded from above and from below near $0$. In particular, for every $N\geqslant2$ and every $\Gamma>0$, $\sup_x|\nabla^2_{x^ix^j}\mathbf E_N(x)|=+\infty$, and the diagonal blocks $\nabla^2_{x^ix^i}\mathbf E_N$ are not bounded from below. This is the same obstruction as for the Bakry-\'Emery criterion.
 
\emph{Otto-Reznikoff.} The criterion of \cite{OttoReznikoff07} requires pointwise bounds $|\nabla^2_{x^ix^j}\mathbf E_N|\le\kappa_{ij}$ on the mixed second derivatives, and the positivity of the matrix built from the one-site log-Sobolev constants and the $\kappa_{ij}$. Here $\kappa_{ij}=+\infty$, and the criterion fails.

\emph{Marton.} Marton's result \cite[Theorem~1]{Marton13} is an approximate tensorization of the relative entropy: under its assumptions, for every probability measure $Q$ on $(\T^3)^N$,
\bbeq\label{eq:AT}
H(Q|\mu_N)\le C\sum_{i=1}^N\int_{(\T^3)^N} H\big(Q_i(\cdot|x^{-i})\,\big|\,\mu_i(\cdot|x^{-i})\big)\,Q(dx),
\nneq
where $Q_i(\cdot|x^{-i})$ is the conditional law of $x^i$ under $Q$. More precisely, \cite[Theorem~1]{Marton13} gives \eqref{eq:AT} with the weights $\rho_i/\rho$, where $\rho_i$ is the log-Sobolev rate of the conditional laws and $\rho$ is the constant of \cite[Assumption~3]{Marton13} (in particular \eqref{eq:AT} holds with $C=\max_i\rho_i/\rho$). Combined with the one-site log-Sobolev inequalities, \eqref{eq:AT} would give a log-Sobolev inequality for $\mu_N$ with constant $C\cdot C_{{\rm LS},1}$, without defect. However, \eqref{eq:AT} requires two assumptions on the Hamiltonian, see \cite[Section~1]{Marton13}. First, the diagonal blocks $\nabla^2_{x^ix^i}\mathbf E_N$ must be bounded from below \cite[Assumption~2]{Marton13}, which fails since they have the eigenvalue $-\frac{\Gamma}{N-1}|x^i-x^k|^{-3}$ near the collisions. Second, the interdependence is measured by the matrix of the pointwise bounds on $\nabla^2_{x^ix^j}\mathbf E_N$, normalized by the one-site log-Sobolev constants, whose norm must be smaller than $1$ \cite[Assumption~3 and (2.3')]{Marton13}, and these bounds are infinite here.

This failure is not a technical artefact. In the proof of Proposition~\ref{P-coef_W2}, the displacement of the pole $x^j$ produces the term $\cov_{\mu_t}(e\cdot\nabla W(\cdot-\eta(t)),f_t)$, in which the singular force is integrated against the conditional law $\mu_t$ of the Gibbs measure itself. The density of $\mu_t$ is bounded and vanishes at the pole. As a result, this term is controlled by the $H^{-1}(\mu_t)$-norm of the force (Lemma~\ref{lem26}) times the energy of $f_t$, and the energy is absorbed by the Hamilton-Jacobi term $-\frac14\mu_t(|\nabla f_t|^2)$. On the contrary, \eqref{eq:AT} must hold for every $Q$. In Marton's proof \cite[Proposition~3, (3.15)]{Marton13}, the displacement of $x^j$ produces a term of the form
$$
\iint_{\T^3\times\T^3}\big(\nabla_{x^i}V(\cdot,\xi)-\nabla_{x^i}V(\cdot,\eta)\big)\,\pi(d\xi,d\eta),
$$
where $\pi$ is a coupling between the conditional laws of $x^j$ under the arbitrary measure $Q$ and under $\mu_N$. Since $Q$ is arbitrary, its conditional laws may concentrate near the poles, and this term can only be controlled through a Lipschitz bound on $\nabla_{x^i}V$ with respect to $x^j$, that is, through the mixed Hessian. In short, the transport argument only tests the interaction against the conditional laws of $\mu_N$, where the singularity is integrable, while the entropy argument tests it against arbitrary laws.

\emph{Zegarlinski.} Following \cite[(47)]{GLWZ2022}, for $i\neq j$, the Zegarlinski interdependence coefficient $c^{\rm Z}_{ij}$ of $\mu_j:=\mu_N(dx^j|x^{-j})$ upon $x^i$ is the best nonnegative constant such that
\bbeq\label{eq:Zeg}
\big|\nabla_{x^i}\big(\mu_j(f^2)\big)^{1/2}\big|\le\big(\mu_j(|\nabla_{x^i}f|^2)\big)^{1/2}+c^{\rm Z}_{ij}\big(\mu_j(|\nabla_{x^j}f|^2)\big)^{1/2}
\nneq
for all smooth positive functions $f$ on $(\T^3)^N$, and Zegarlinski's criterion \cite[Theorem~0.1]{Zegarlinski92} (see \cite[Theorem~16]{GLWZ2022}) requires $\sup_i\max\{\sum_jc^{\rm Z}_{ji},\sum_jc^{\rm Z}_{ij}\}<1$. These coefficients do not involve the Hessian of the Hamiltonian, but the gradient dependence of the conditional laws. For smooth mean-field interactions, $c^{\rm Z}_{ij}=O(\frac1N)$ \cite[Lemma~18]{GLWZ2022}. For the Coulomb gas, this is not the case, as we now show.

Let $f=f(x^j)\geqslant0$ depend only on $x^j$. Then $\nabla_{x^i}f=0$, and \eqref{eq:Zeg} (applied to $f+\varepsilon$, with $\varepsilon\to0$) reduces to
$$
\big|\nabla_{x^i}\mu_j(f^2)\big|\le2\,c^{\rm Z}_{ij}\big(\mu_j(f^2)\,\mu_j(|\nabla_{x^j}f|^2)\big)^{1/2}.
$$ 
Set $c:= {\Gamma}/{(N-1)}$. Since only the factor $e^{-cW(x^j-x^i)}$ of the density of $\mu_j$ depends on $x^i$, we have, for every $f$ depending only on $x^j$, $\nabla_{x^i}\mu_j(f^2)=c\,\cov_{\mu_j}\big(f^2,\nabla W(\cdot-x^i)\big)$. Fix $a\in\T^3$ and a unit vector $e\in\rr^3$, and take
$$
f(x^j):=\varphi\big((x^j-a)/c\big),
$$
where $\varphi\geqslant0$, $\varphi\not\equiv0$, is a fixed smooth function supported in $\{y:\ 1<|y|<2,\ e\cdot y>\frac12\}$. This function only depends on $x^j$. We evaluate the reduced inequality at a configuration $x^{-j}$ such that $x^i=a$ and the other particles $x^k$, $k\neq i,j$, are at distance at least $\frac14$ from $a$. The function $f$ is supported in
$$
A_c:=\Big\{x^j:\ c<|x^j-a|<2c,\ e\cdot(x^j-a)>\tfrac c2\Big\},\quad\text{and}\quad |\nabla f|\le C_0\,c^{-1},\quad C_0:=\|\nabla\varphi\|_\infty .
$$
Since $W(y)=|y|^{-1}+R(y)$, we have $cW(x^j-a)\in(\frac12,1)+O(c)$ on $A_c$, and the contribution $c\sum_{k\neq i,j}W(x^j-x^k)$ of the other particles varies by $O(\Gamma c)$ on $A_c$. Thus the density of $\mu_j$ is comparable to a constant on $A_c$, uniformly in $N$, and $\mu_j(A_c)\le K\mu_j(f^2)$ for some $K$ depending only on $\varphi$. Moreover, with $y=x^j-a$, we have $e\cdot y>\frac c2$ and $|y|<2c$ on $A_c$, so that
$$
e\cdot\nabla W(y)=-\frac{e\cdot y}{|y|^3}+O(1)\le-\frac{1}{16}\,c^{-2}+O(1).
$$
Using $|\nabla_{x^i}\mu_j(f^2)|\geqslant|e\cdot\nabla_{x^i}\mu_j(f^2)|$ and $|\mu_j(e\cdot\nabla W(\cdot-a))|\le e^{-\Gamma W_*}\|\nabla W\|_{L^1(\T^3,dx)}$, we get that, at $x^i=a$,
$$
\big|\nabla_{x^i}\mu_j(f^2)\big|\geqslant c\,\Big|\mu_j\big(f^2\,e\cdot\nabla W(\cdot-a)\big)-\mu_j(f^2)\,\mu_j\big(e\cdot\nabla W(\cdot-a)\big)\Big|\geqslant\Big(\frac{1}{16}-Cc^2\Big)c^{-1}\mu_j(f^2),
$$
for some $C$ depending only on $W$ and $\Gamma$, whereas, since $|\nabla f|^2\le C_0^2c^{-2}\,\mathbf 1_{A_c}$,
$$
\big(\mu_j(f^2)\,\mu_j(|\nabla_{x^j}f|^2)\big)^{1/2}\le C_0\big(\mu_j(f^2)\,\mu_j(A_c)\big)^{1/2}c^{-1}\le C_0K^{1/2}\,c^{-1}\mu_j(f^2).
$$
Plugging these bounds into the reduced inequality, we deduce that $c^{\rm Z}_{ij}\geqslant\kappa_0:=\frac{1}{64\,C_0K^{1/2}}$ as soon as $c$ is small enough, i.e.\ for $N\geqslant N_0(\Gamma)$, where $\kappa_0$ does not depend on $N$ nor on $\Gamma$. Since this holds for every pair $i\neq j$, $\sum_jc^{\rm Z}_{ij}\geqslant(N-1)\kappa_0\to\infty$, and Zegarlinski's condition fails for $N$ large, however small $\Gamma$ is.

\emph{The route we follow.} We therefore follow a different route, suggested by Rosenzweig \cite[Section~5]{Rosenzweig_partI}. We do not compare the conditional laws of $\mu_N$ with each other. Instead, we compare $\mu_N$ with the free product measure $\mathbf m _N$, for which the entropy satisfies the independent tensorization, and the interaction only appears through an energy term. Its exponential moments are controlled by a uniform bound on the partition function. The price to pay is an additive defect in the log-Sobolev inequality, which is removed by Rothaus' lemma and the uniform Poincar\'e inequality of Theorem~\ref{T_W2I}\emph{(1)}.

\subsection{One-site log-Sobolev inequality}\label{sec:LSI1}

\bprop\label{LSI_single}
Let $\mu^{\y}$ be the one-particle measure \eqref{mu_single} on $\T^3$, with total charge $\bar\gamma<2\pi$. Then, for every $f\in H^1(\mu^{\y})$,
\bbeq\label{LSI_single1}
\Ent_{\mu^{\y}}(f^2)\le 2C_{{\rm LS},1}\int_{\T^3}|\nabla f|^2\,d\mu^{\y},\qquad C_{{\rm LS},1}:=\frac12\big(A+(B+2)\,C_{{\rm P},1}\big),
\nneq
where $C_{{\rm P},1}=\frac1{2\pi(2\pi-\bar\gamma)}$ is the constant of Proposition~\ref{PI_single} and
$$
A:=\frac32+\bar\gamma\,\|W\|_{L^{2}(\T^3,dx)}C_{{\rm So}},\qquad
B:=2\pi\bar\gamma A+\bar\gamma\,\|W\|_{L^{2}(\T^3,dx)}C_{{\rm So}}+\bar\gamma|W_*|+\frac32\log C_{{\rm So}} .
$$
The constant $C_{{\rm LS},1}$ depends only on $\bar\gamma$, and not on $N$, on the number of poles or on their positions.
\nprop

Recall that $\|W\|_{L^2(\T^3,dx)}\approx1.29$ and $|W_*|\approx0.802$ (see Section~\ref{sec.N}), and that $C_{{\rm So}}\le\frac73$ (see Lemma~\ref{le.Sobolev}). Since the singular density of $\mu^{\y}$ is not bounded from below, the bounded perturbation principle of Holley and Stroock does not apply. The proof combines the ground state identity of Proposition~\ref{PI_single} with the Sobolev inequality \eqref{SobolevI}: it first gives a defective log-Sobolev inequality, and the defect is then removed by the one-particle Poincar\'e inequality. This is, at the level of one particle, the same structure as the proof of Theorem~\ref{T_LSI_intro}.

\bprf
Write $\mu=\mu^{\y}$, $V=V^{\y}$ and $\varphi=Z^{-1/2}e^{-V/2}$. Let $f\in H^1(\mu)$ with $\int_{\T^3} f^2\,d\mu=1$ (the inequality is $2$-homogeneous), and set $g:=f\varphi$. Then, $\|g\|_2=1$ and, by \eqref{eq:H1},
\bbeq\label{eq:gradg}
\|\nabla g\|_2^2\le\EE_\mu[f]+2\pi\bar\gamma .
\nneq

\textbf{Step 1. Defective inequality.} Since $f^2=g^2/\varphi^2$ and $-\log\varphi^2=V+\log Z$, we have
$$
\Ent_\mu(f^2)=\int_{\T^3}g^2\log g^2\,dx+\int_{\T^3}Vg^2\,dx+\log Z .
$$
All terms are finite, since $g\in L^6$ by \eqref{SobolevI} and $V\in L^{3/2}$. By Jensen's inequality for the concave function $\log$ under the probability measure $g^2dx$, by \eqref{SobolevI} and $\|g\|_2=1$, and by $\log(1+u)\le u$, we get that:
$$
\aligned
\int_{\T^3}g^2\log g^2\,dx&=\frac12\int_{\T^3}\log(g^4)\,g^2dx\le\frac12\log\|g\|_6^6=\frac32\log\|g\|_6^2\\
&\le\frac32\log C_{{\rm So}}+\frac32\log\big(1+\|\nabla g\|_2^2\big)\le\frac32\log C_{{\rm So}}+\frac32\|\nabla g\|_2^2.
\endaligned
$$
By H\"older's inequality, $\|V\|_{3/2}\le\bar\gamma\|W\|_{3/2}\le\bar\gamma\|W\|_2$ (since $dx$ is a probability measure), and by \eqref{SobolevI},
$$
\int_{\T^3} Vg^2\,dx\le\|V\|_{3/2}\|g\|_6^2\le\bar\gamma\|W\|_2C_{{\rm So}}\big(\|\nabla g\|_2^2+1\big).
$$
Finally $\log Z\le\bar\gamma|W_*|$ by \eqref{eq:Zbound}. Summing these three bounds, we obtain
$$
\Ent_\mu(f^2)\le\Big(\frac32+\bar\gamma\|W\|_2C_{{\rm So}}\Big)\|\nabla g\|_2^2+\bar\gamma\|W\|_2C_{{\rm So}}+\bar\gamma|W_*|+\frac32\log C_{{\rm So}},
$$
and \eqref{eq:gradg} gives $\Ent_\mu(f^2)\le A\,\EE_\mu[f]+B$. By homogeneity, for every $f\in H^1(\mu)$,
\bbeq\label{defLSI_single}
\Ent_\mu(f^2)\le A\,\EE_\mu[f]+B\int_{\T^3}f^2\,d\mu .
\nneq

\textbf{Step 2. Rothaus' lemma.} By Rothaus' lemma \cite{Rothaus85} (see also \cite[(5.29)]{Rosenzweig_partI}), with $\tilde f:=f-\mu(f)$,
\bbeq\label{eq:Rothaus}
\Ent_\mu(f^2)\le\Ent_\mu(\tilde f^2)+2\var_\mu(f).
\nneq
Applying \eqref{defLSI_single} to $\tilde f$ (note that $\EE_\mu[\tilde f]=\EE_\mu[f]$ and $\int_{\T^3}\tilde f^2\,d\mu=\var_\mu(f)$), and then Proposition~\ref{PI_single}, we deduce that:
$$
\Ent_\mu(f^2)\le A\,\EE_\mu[f]+(B+2)\var_\mu(f)\le\big(A+(B+2)C_{{\rm P},1}\big)\EE_\mu[f]=2C_{{\rm LS},1}\EE_\mu[f].
$$
The proof of the proposition is complete. 
\nprf


A compactness argument shows that a uniform one-site log-Sobolev inequality also holds for any bounded total charge, with a non-explicit constant; we shall not need it, since the propagation step is perturbative anyway.

\subsection{An abstract propagation principle with a defect}\label{sec:defect}

The following statement isolates the assumptions needed to propagate the one-site log-Sobolev inequality to $\mu_N$. It holds for any even kernel $W$ on $\T^d$ for which the framework of Section~\ref{subsec:FI} applies. For simplicity, we state it on $\T^3$.

\bthm\label{T-abstract-LSI}
Assume that $\int_{\T^3}W\,dx=0$ and that $W$ is bounded from below. Let $\mu_N$ be given by \eqref{Coulomb-measure}, and assume that:
\begin{enumerate}
  \item[{\rm (\textbf{H1})}] {\rm (One-site LSI)} There is $C_{{\rm LS},1}<\infty$ such that, for all $N$, $i$ and $x^{-i}$, $\Ent_{\mu_i(\cdot|x^{-i})}(g^2)\le2C_{{\rm LS},1}\int_{\T^3}|\nabla g|^2\,d\mu_i(\cdot|x^{-i})$ for every $g\in H^1(\mu_i(\cdot|x^{-i}))$;
  \item[{\rm (\textbf{H2})}] {\rm (Partition bound)} There is $\lambda>1$ such that $\Lambda_\lambda:=\sup_{N\geqslant2}\log Z_N(1+\lambda)<\infty$ (see \eqref{eq:energy}).
\end{enumerate}
Then, for every $N\geqslant2$ and every $f\in H^1(\mu_N)$, the following defective log-Sobolev inequality holds:
\bbeq\label{eq:defLSI_abstract}
\Ent_{\mu_N}(f^2)\le\frac{\lambda}{\lambda-1}\,2C_{{\rm LS},1}\int_{(\T^3)^N}|\nabla f|^2\,d\mu_N+D_\lambda\int_{(\T^3)^N} f^2\,d\mu_N,\qquad D_\lambda:=\frac{2\lambda}{\lambda^2-1}\,\Lambda_\lambda .
\nneq
If moreover
\begin{enumerate}
  \item[{\rm (\textbf{H3})}] {\rm (Uniform Poincar\'e inequality)} $\var_{\mu_N}(f)\le C_P\int_{(\T^3)^N}|\nabla f|^2\,d\mu_N$ for all $N$ and $f\in H^1(\mu_N)$,
\end{enumerate}
then $\mu_N$ satisfies the log-Sobolev inequality $\Ent_{\mu_N}(f^2)\le 2C_{\rm LS}\int_{(\T^3)^N}|\nabla f|^2\,d\mu_N$, with
\bbeq\label{eq:CLS_abstract}
2C_{\rm LS}=\frac{\lambda}{\lambda-1}\,2C_{{\rm LS},1}+(D_\lambda+2)\,C_P .
\nneq
\nthm

Let us comment on the assumptions. No interdependence coefficient appears in (\textbf{H1})-(\textbf{H2}). The interaction enters only through the energy $\mathbf E_N$ and its exponential moments under the free product measure (\textbf{H2}). Interdependence coefficients are still used, but only for the variance in (\textbf{H3}), through Dobrushin's coefficients in total variation, which require nothing more than $W\in L^1$ (Lemma~\ref{DIC}). The defect $D_\lambda$ comes from the interaction energy; for a product measure ($W=0$), $\Lambda_\lambda=0$ for every $\lambda>1$, and letting $\lambda\to\infty$ in \eqref{eq:defLSI_abstract} gives back the independent tensorization. The factor $\frac\lambda{\lambda-1}$ comes from the use of the Donsker-Varadhan formula below.

The proof follows \cite[Section~5]{Rosenzweig_partI}. In Rosenzweig's planar setting, the Gibbs measure is first recentered at the thermal equilibrium. On $\T^3$, since $\int_{\T^3}W\,dx=0$, no recentering is needed, and the reference measure is simply $\mathbf m_N$.


\bprf
\textbf{Step 0. Reduction.} By Proposition~\ref{prop:Wu}\emph{(c)}, Fatou's lemma for the entropy (since $u\log u\geqslant-e^{-1}$), and the convergence of the right-hand side, it is enough to consider $f\in C_c^\infty(\mathcal D)$, where $\mathcal D=(\T^3)^N\setminus\SS_N$. By homogeneity, we assume that $\int_{(\T^3)^N} f^2\,d\mu_N=1$, and we set 
$$Q:=f^2\mu_N.$$ Note that $H(Q|\mu_N)=\Ent_{\mu_N}(f^2)$. All the quantities below are finite. Indeed, since $f$ is bounded and $\mathbf E_N\geqslant\frac{\Gamma N}{2}W_*$, the density of $Q$ with respect to $\mathbf m _N$ is bounded by $\|f\|_\infty^2e^{-\Gamma NW_*/2}$ (recall $Z_N(1)\geqslant1$), and $\mathbf E_N\in L^1(\mathbf m _N)$ since $W\in L^1(\T^3,dx)$. For the same reason, for $Q$-a.e.\ $x$, the conditional laws $Q_i(\cdot|x^{-i})$ defined below have bounded densities with respect to $dz$ (throughout the proof, $z\in\T^3$ denotes the variable $x^i$, and $dz$ the Lebesgue measure on $\T^3$). For $i\in\{1,\dots,N\}$, let $Q_i(\cdot|x^{-i})$ be the conditional law of $x^i$ under $Q$, and set
$$
S_N(Q):=\sum_{i=1}^N\int_{(\T^3)^N} H\big(Q_i(\cdot|x^{-i})\,\big|\,\mu_i(\cdot|x^{-i})\big)\,Q(dx).
$$

The quantity $S_N(Q)$ is the sum of the averaged conditional entropies of $Q$ with respect to the one-site measures. It is the bridge between the global entropy and the one-site inequalities: we bound it from below by $H(Q|\mu_N)$, up to the interaction term $\ee_Q[\mathbf E_N]$ (Steps~1 and~2), which gives a defective decomposition inequality for the relative entropy, and from above by the one-site log-Sobolev inequalities (Step~3). Comparing the two bounds gives the defective log-Sobolev inequality \eqref{eq:defLSI_abstract}.

\smallskip

\textbf{Step 1. Independent tensorization with respect to the free measure.} The conditional law $\mu_i(dz|x^{-i})$ has density $e^{-\Phi_i(z)}/\int_{\T^3} e^{-\Phi_i}\,dz$ with respect to $dz$, where $\Phi_i(z):=\frac{\Gamma}{N-1}\sum_{j\neq i}W(z-x^j)$. Since $\int_{\T^3}\Phi_i\,dz=0$, Jensen's inequality gives $\int_{\T^3} e^{-\Phi_i}\,dz\geqslant1$, which yields
$$
\aligned
H\big(Q_i(\cdot|x^{-i})\,|\,\mu_i(\cdot|x^{-i})\big)&=H\big(Q_i(\cdot|x^{-i})\,|\,dz\big)+\int_{\T^3}\Phi_i\,dQ_i(\cdot|x^{-i})+\log\int_{\T^3} e^{-\Phi_i}\,dz\\
&\geqslant H\big(Q_i(\cdot|x^{-i})\,|\,dz\big)+\int_{\T^3}\Phi_i\,dQ_i(\cdot|x^{-i}).
\endaligned
$$
We integrate with respect to $Q(dx)$ and sum over $i$. The relative entropy with respect to the product measure $\mathbf m _N$ is subadditive, i.e.  
\bbeq\label{eq:subadd}
H(Q|\mathbf m _N)\le\sum_{i=1}^N\int_{(\T^3)^N} H(Q_i(\cdot|x^{-i})|dz)\,Q(dx).
\nneq 
Moreover, $\sum_{i=1}^N\Phi_i(x^i)=2\mathbf E_N(x)$. Consequently,
 $
S_N(Q)\geqslant H(Q|\mathbf m _N)+2\,\ee_Q[\mathbf E_N]$. 
Since $d\mu_N/d\mathbf m _N=e^{-\mathbf E_N}/Z_N(1)$, we have $H(Q|\mathbf m _N)=H(Q|\mu_N)-\ee_Q[\mathbf E_N]-\log Z_N(1)$, and we get that:
\bbeq\label{eq:SNlow}
S_N(Q)\geqslant H(Q|\mu_N)+\ee_Q[\mathbf E_N]-\log Z_N(1).
\nneq

\textbf{Step 2. The interaction term and the defective decomposition inequality.} Let $\lambda>1$ be given by (\textbf{H2}). 
The energy $\ee_Q[\mathbf E_N]$ may be negative, and we bound it from below by the Donsker-Varadhan variational formula $\ee_Q[F]\le H(Q|\mu_N)+\log\int_{(\T^3)^N} e^{F}\,d\mu_N$, valid for every measurable function $F$ bounded from above, applied to $F=-\lambda \mathbf E_N$ (recall that $\mathbf E_N\geqslant\frac{\Gamma N}2W_*$ over $(\T^3)^N$,   $\ee_Q[|\mathbf E_N|]<\infty$ by Step $0$). Since $\int_{(\T^3)^N} e^{-\lambda \mathbf E_N}\,d\mu_N=Z_N(1+\lambda)/Z_N(1)$, we obtain
$$
\ee_Q[\mathbf E_N]\geqslant-\frac1\lambda H(Q|\mu_N)-\frac1\lambda\log Z_N(1+\lambda)+\frac1\lambda\log Z_N(1).
$$
Plugging this into \eqref{eq:SNlow}, we deduce that:
$$
S_N(Q)\geqslant\Big(1-\frac1\lambda\Big)H(Q|\mu_N)-\frac1\lambda\log Z_N(1+\lambda)-\Big(1-\frac1\lambda\Big)\log Z_N(1).
$$
By Jensen's inequality for the convex function $u\mapsto u^{1+\lambda}$ under the probability measure $\mathbf m _N$, $Z_N(1)\le Z_N(1+\lambda)^{1/(1+\lambda)}$, and since $\log Z_N(1)\geqslant0$, we get that $(1-\frac1\lambda)\log Z_N(1)\le\frac{\lambda-1}{\lambda(\lambda+1)}\log Z_N(1+\lambda)$. Multiplying by $\frac{\lambda}{\lambda-1}$, we obtain (recall (\textbf{H2})):
\bbeq\label{eq:defdecomp}
H(Q|\mu_N)\le\frac{\lambda}{\lambda-1}S_N(Q)+\Big(\frac1{\lambda-1}+\frac1{\lambda+1}\Big)\log Z_N(1+\lambda)\le\frac{\lambda}{\lambda-1}S_N(Q)+D_\lambda .
\nneq
This is the defective decomposition inequality for the relative entropy: up to the factor $\frac{\lambda}{\lambda-1}$ and the additive defect $D_\lambda$, the global entropy $H(Q|\mu_N)$ is bounded by the sum $S_N(Q)$ of the averaged conditional entropies.

\textbf{Step 3. Upper bound on $S_N(Q)$ and defective log-Sobolev inequality.} For $\mu_N$-a.e.\ $x^{-i}$, $$Q_i(dx^i|x^{-i})=\frac{f^2(x^i,x^{-i})}{\mu_i(f^2|x^{-i})}\, \mu_i(dx^i|x^{-i}) \text{ i.e. $
\frac{dQ_i(\cdot| x^{-i})}
{d\mu_i(\cdot| x^{-i})}(x^i)
=
\frac{f^2(x^i,x^{-i})}
{\mu_i(f^2| x^{-i})}$},$$
 and the law $Q^{-i}(dx^{-i})$ of $x^{-i}$ under $Q$ has density $\mu_i(f^2|x^{-i})$ with respect to its law $\mu_N^{-i}$ under $\mu_N$ ($\mu_N^{-i}(dx^{-i})
=
\int_{\mathbb{T}^3}
\mu_N(dx^i,dx^{-i})$ and $\mu_i(f^2|x^{-i})
=
\int_{\mathbb{T}^3}
f^2(x^i,x^{-i})\,\mu_i(dx^i| x^{-i})$). Indeed, we have by the disintegration $\mu_N(dx^i,dx^{-i})=\mu_i(dx^i|x^{-i})\,\mu_N^{-i}(dx^{-i})$, we have   
 $$Q(dx^i,dx^{-i})=\frac{f^2\,\mu_i(dx^i|x^{-i})}{\mu_i(f^2|x^{-i})}\,\mu_i(f^2|x^{-i})\,\mu_N^{-i}(dx^{-i})= Q_i(dx^i|x^{-i})\, \mu_i(f^2|x^{-i})\,\mu_N^{-i}(dx^{-i}).$$
 Thus, we get that for all $x^{-i}$:
$$ \begin{aligned}
H\bigl(Q_i(\cdot| x^{-i})
|\mu_i(\cdot| x^{-i})\bigr)&=
\int_{\mathbb T^3}
\log\left(
\frac{f^2(x^i,x^{-i})}
{\mu_i(f^2| x^{-i})}
\right)
Q_i(dx^i| x^{-i})\\
&=
\frac{1}{\mu_i(f^2| x^{-i})}
\int_{\mathbb T^3}
f^2(x^i,x^{-i})
\log\left(
\frac{f^2(x^i,x^{-i})}
{\mu_i(f^2| x^{-i})}
\right)
\mu_i(dx^i| x^{-i})\\
&= \frac{\Ent_{\mu_i(\cdot|x^{-i})}(f^2)}{\mu_i(f^2| x^{-i})}.
\end{aligned}
$$
Because, $Q^{-i}(dx^{-i}) =\mu_i(f^2|x^{-i})\,\mu_N^{-i}(dx^{-i})$ we then deduce that 
\begin{align*}
\int_{(\T^3)^N} H\big(Q_i(\cdot|x^{-i})\, |\,\mu_i(\cdot|x^{-i})\big)\,Q(dx)&= \int_{(\T^3)^{N-1}} H\big(Q_i(\cdot|x^{-i})\, |\,\mu_i(\cdot|x^{-i})\big)\,Q^{-i}(dx^{-i})\\
&= \int_{(\T^3)^N}\Ent_{\mu_i(\cdot|x^{-i})}(f^2)\,d\mu_N(x).
\end{align*}
  Moreover, for every $x^{-i}$, the function $z\mapsto f(z,x^{-i})$ is smooth on $\T^3$ and vanishes in a neighborhood of the poles $\{x^j,\ j\neq i\}$, since the support of $f$ is a compact subset of $\mathcal D$. It thus belongs to $H^1(\mu_i(\cdot|x^{-i}))$. Therefore, by (\textbf{H1}) applied to $z\mapsto f(z,x^{-i})$,
$$
 \sum_{i=1}^N\int_{(\T^3)^N}\Ent_{\mu_i(\cdot|x^{-i})}(f^2)\,d\mu_N(x)\le 2C_{{\rm LS},1}\sum_{i=1}^N\int_{(\T^3)^N}|\nabla_{i}f|^2\,d\mu_N=2C_{{\rm LS},1}\int_{(\T^3)^N}|\nabla f|^2\,d\mu_N .
$$
Therefore, $S_N(Q)\le 2C_{{\rm LS},1}\int_{(\T^3)^N}|\nabla f|^2\,d\mu_N$. Plugging this bound into \eqref{eq:defdecomp}, and since $H(Q|\mu_N)=\Ent_{\mu_N}(f^2)$ and $\int_{(\T^3)^N}f^2\,d\mu_N=1$, we arrive at the defective log-Sobolev inequality \eqref{eq:defLSI_abstract}.

\smallskip

\textbf{Step 4. Conclusion.} Under (\textbf{H3}), Rothaus' lemma \eqref{eq:Rothaus} for $\mu_N$, \eqref{eq:defLSI_abstract} applied to $\tilde f=f-\mu_N(f)$ and (\textbf{H3}) give
$
\Ent_{\mu_N}(f^2)\le\frac{\lambda}{\lambda-1}2C_{{\rm LS},1}\int_{(\T^3)^N}|\nabla f|^2\,d\mu_N+(D_\lambda+2)\var_{\mu_N}(f)\le2C_{\rm LS}\int_{(\T^3)^N}|\nabla f|^2\,d\mu_N$, the desired result. 
\nprf

\subsection{The partition function bound}\label{sec:partition}
The key input to apply Theorem~\ref{T-abstract-LSI} is the partition bound (\textbf{H2}). For the Coulomb gas on $\T^3$, it holds at every temperature, by a result of Wang and Zhao \cite{WangZhao26}.

\blem\label{le.partition} {\rm (Wang and Zhao \cite[Theorem~1.1]{WangZhao26})}
For every $t>0$,
\bbeq\label{eq:partition}
\Lambda(t):=\sup_{N\geqslant2}\log Z_N(t)<\infty .
\nneq
Moreover, $\Lambda(t)\to0$ as $\Gamma\to0$, for every fixed $t>0$ (see \cite[Theorem~1.1]{WangZhao26}).
\nlem

The result of \cite[Theorem~1.1]{WangZhao26} is stated on $\rr^d$, but, as indicated to us by X.~Zhao, its proof applies on $\T^d$ with a simple adaptation, for the periodic logarithmic and Riesz kernels in the range $0\le s<\frac d2$, which includes the Coulomb kernel on $\T^2$ and $\T^3$.

The condition $W\in L^2$ is sharp: if $W\notin L^2(\T^d)$, then $Z_N(t)\to\infty$ for every $t>0$ \cite[Proposition~2.2(ii)]{DJ26}. At high temperature, the bound \eqref{eq:partition} was obtained before by Duerinckx and Jabin \cite[Theorem~2.1(i)]{DJ26}, and independently by Rosenzweig (see \cite[Section~2]{DJ26}). For repulsive logarithmic-type interactions, a uniform bound at every temperature, for a bounded reference density, is due to Delgadino and Gvalani \cite[Theorem~2.5]{DG25}.

\subsection{Proof of Theorem~\ref{T_LSI_intro}}\label{sec:LSI_proof}

\bthm\label{T_LSI}
Assume that $r_{\rm D}<1$, i.e.\ \eqref{T_W2I1} holds. Then, for every $N\geqslant2$ and every $f\in H^1(\mu_N)$,
\bbeq\label{eq:LSI}
\Ent_{\mu_N}(f^2)\le 2C_{\rm LS}\int_{(\T^3)^N}|\nabla f|^2\,d\mu_N,\qquad C_{\rm LS}:=2C_{{\rm LS},1}+\frac{D+2}{2(1-r_{\rm D})\cdot2\pi(2\pi-\Gamma)},
\nneq
where $C_{{\rm LS},1}$ is given by Proposition~\ref{LSI_single} with $\bar\gamma=\Gamma$, and $D:=\frac43\Lambda(3)<\infty$, with $\Lambda(3)=\sup_{N\geqslant2}\log Z_N(3)$ given by Lemma~\ref{le.partition}.
\nthm

\bprf
Recall that $r_{\rm D}<1$ forces $\Gamma<\pi$. We apply Theorem~\ref{T-abstract-LSI} with $\lambda=2$. Assumption (\textbf{H1}) holds with the constant $C_{{\rm LS},1}$ of Proposition~\ref{LSI_single}, since each $\mu_i(\cdot|x^{-i})$ has total charge $\Gamma<2\pi$. Assumption (\textbf{H2}) holds with $\lambda=2$ and $\Lambda_2=\Lambda(3)<\infty$ by Lemma~\ref{le.partition}, which gives $D_2=\frac43\Lambda_2=D$. Assumption (\textbf{H3}) holds with $C_P=\frac{1}{(1-r_{\rm D})2\pi(2\pi-\Gamma)}$ by Theorem~\ref{T_W2I}\emph{(1)}. Then, \eqref{eq:CLS_abstract} gives $2C_{\rm LS}=4C_{{\rm LS},1}+(D+2)C_P$, which is \eqref{eq:LSI}.
\nprf


%

\subsection{Extensions and discussion}\label{sec:LSI_ext}

\subsubsection{Extensions}
Theorem~\ref{T-abstract-LSI} applies to the two-dimensional Coulomb gas and to periodic Riesz gases.

\bcor\label{cor:LSI_ext}
Concerning the two-dimensional Coulomb gas and to periodic Riesz gases, we have the following: 
\begin{enumerate}[(a)]
  \item {\rm (Two-dimensional Coulomb gas)} Let $W$ be given by \eqref{eq:W2D}, and assume that $r^{\scriptscriptstyle (2D)}_{\rm D}<1$ (see Theorem \ref{T_2D}\emph{(1)}). Then,   $\mu_N$ satisfies a (uniform in $N$) log-Sobolev inequality on $(\T^2)^N$ with constant
  $$
  C_{\rm LS}=2C^{\scriptscriptstyle (2D)}_{{\rm LS},1}+\frac{D+2}{2(1-r^{\scriptscriptstyle (2D)}_{\rm D})\cdot\pi(4\pi-\Gamma)},
  $$
  where $D:=\frac43\sup_{N\geqslant2}\log Z_N(3)<\infty$, $C^{\scriptscriptstyle (2D)}_{{\rm LS},1}=\frac12(A+(B+2)C^{\scriptscriptstyle (2D)}_{{\rm P},1})$, with $C^{\scriptscriptstyle (2D)}_{{\rm P},1}=\frac1{\pi(4\pi-\Gamma)}$, $A=2+\Gamma\|W\|_{L^2(\T^2,dx)}C^{\scriptscriptstyle (2D)}_{{\rm So}}$ and $B=\pi\Gamma A+\Gamma\|W\|_{L^2(\T^2,dx)}C^{\scriptscriptstyle (2D)}_{{\rm So}}+\Gamma|W_*|+2\log C^{\scriptscriptstyle (2D)}_{{\rm So}}$.

  \item {\rm (Periodic Riesz gases)} Let $0<s\le1$, and assume that $r_{\rm D}^{(s)}<1$ (see Theorem \ref{T_Riesz}\emph{(1)}), $\Gamma<8\pi^2/K_s$ (where we recall that $K_s=\sup_{\T^3\setminus\{0\}}\Delta W_s$). Then, $\mu_N^{\scriptscriptstyle (s)}$ satisfies a (uniform in $N$) log-Sobolev inequality on $(\T^3)^N$ with constant
  $$
  C_{\rm LS}=2C^{(s)}_{{\rm LS},1}+\frac{(D+2)\,C^{(s)}_{{\rm P},1}}{2(1-r_{\rm D}^{(s)})},
  $$
  where $D:=\frac43\sup_{N\geqslant2}\log Z_N(3)<\infty$, $C^{(s)}_{{\rm LS},1}=\frac12(A+(B+2)C^{(s)}_{{\rm P},1})$, with $A=\frac32+\Gamma\|W_s\|_{L^2}C_{{\rm So}}$ and $B=\frac12\Gamma K_sA+\Gamma\|W_s\|_{L^2}C_{{\rm So}}+\Gamma|W_{s,*}|+\frac32\log C_{{\rm So}}$.
\end{enumerate}
\ncor

\bprf
We apply Theorem~\ref{T-abstract-LSI} with $\lambda=2$. Assumption (\textbf{H3}) is given by Theorems~\ref{T_2D}\emph{(1)} and~\ref{T_Riesz}\emph{(1)}, and (\textbf{H2}) by the remark following Lemma~\ref{le.partition} (see \cite{WangZhao26}), since $W\in L^2(\T^2,dx)$ and $W_s\in L^2(\T^3,dx)$ for $s<\frac32$. For (\textbf{H1}), the proof of Proposition~\ref{LSI_single} applies with the following changes. The bound \eqref{eq:gradg} becomes $\|\nabla g\|_2^2\le\EE_{\mu^{\y}}[f]+\frac{\kappa\bar\gamma}2$ (for $\int f^2\,d\mu^{\y}=1$), with $\kappa=2\pi$ in case~(a) and $\kappa=K_s$ in case~\emph{(b)}. In case~\emph{(a)}, we use \eqref{Sobolev2D} instead of \eqref{SobolevI}: Jensen's inequality (with the probability measure $
\nu(dx)=g^2(x)\,dx$), \eqref{Sobolev2D} and $\log(1+u)\le u$ give $\int_{\T^2} g^2\log g^2\,dx\le2\log\|g\|_4^2\le2\log C^{\scriptscriptstyle (2D)}_{{\rm So}}+2\|\nabla g\|_2^2$, H\"older's inequality gives $\int_{\T^2}|V|g^2\,dx\le\|V\|_2\|g\|_4^2\le\Gamma\|W\|_2C^{\scriptscriptstyle (2D)}_{{\rm So}}(\|\nabla g\|_2^2+1)$, and $\log Z\le\Gamma|W_*|$. Together with $\|\nabla g\|_2^2\le\EE_{\mu^{\y}}[f]+\pi\Gamma$, this gives the constants $A$ and $B$ of case~\emph{(a)}.
 Case~~\emph{(b)} is identical to Proposition~\ref{LSI_single}, with $\|W_s\|_{3/2}\le\|W_s\|_2$ (recall also that    $\|\nabla g\|_2^2\le\EE_\mu[f]+\frac12\int_{\T^3}g^2\,\Delta V $, and thus, since $\Delta V\le\Gamma K_s$ off the poles,   the last term is at most $\frac12\Gamma K_s$, which explains the first term in $B$). 
\nprf

\subsubsection{Dependent tensorization versus tensorization with a defect}
\label{sec.discu}
Let us compare the two propagation principles and the assumptions on the interaction they require.
\smallskip

\emph{Dependent tensorization} (Theorem~\ref{T-abstract-W}) propagates a one-site inequality to $\mu_N$ through interdependence coefficients, i.e.\ through the Lipschitz dependence of the conditional laws on the other coordinates, measured by a transport distance. The global constant is the one-site constant divided by $1-r$. Since $r\to0$ as $\Gamma\to0$, it tends to the one-site constant: for the Poincar\'e inequality, the constant $\frac{1}{(1-r_{\rm D})\cdot2\pi(2\pi-\Gamma)}$ of Theorem~\ref{T_W2I}\emph{(1)} tends to the optimal constant $1/({4\pi^2})$ of the Lebesgue measure on $\T^3$, and for the $W_2I$ inequality, the constant $C_{W_2I}= {C_{W_2I,1}}/({1-r_{W_2}})$ of Theorem~\ref{T_W2I}\emph{(2)} tends to the one-site constant $C_{W_2I,1}$ at $\Gamma=0$. In this sense, the method is sharp in the weak-coupling limit (when $\Gamma\to0$, the interaction disappears). 

 For the Poincar\'e inequality, Dobrushin's coefficients require $W\in L^1$. For the $W_2I$ inequality, the $W_2$-interdependence coefficients require a uniform bound on the $H^{-1}(\mu^{\y})$-norm of $\nabla W$, which holds when $\nabla W\in L^{2d/(d+2)}(\T^d)$ ($d\geqslant3$). All the constants are explicit. For the log-Sobolev inequality, the available criteria (Zegarlinski, Otto-Reznikoff, Marton) fail for singular kernels, however small $\Gamma$ is (Section \ref{sec:obstruction}). 
\smallskip

\emph{Tensorization with a defect} (Theorem~\ref{T-abstract-LSI}) does not compare the conditional laws with each other. It uses the independent tensorization of the entropy with respect to the free product measure (see \eqref{eq:subadd} in Step~1 of the proof of Theorem~\ref{T-abstract-LSI}), and the interaction only enters through the exponential moments of the energy, i.e.\ the partition bound (\textbf{H2}), which holds at every temperature when $W\in L^2$ (Lemma~\ref{le.partition}). The price is an additive defect $D_\lambda$ (see \eqref{eq:defLSI_abstract}), which requires a uniform Poincar\'e inequality (\textbf{H3}) to be removed (here Theorem~\ref{T_W2I}\emph{(1)}). In a addition, a global constant which does not tend to the one-site constant in the limit  $\Gamma\to0$ remains. Indeed, as $\Gamma\to0$, $D_\lambda\to0$ and $C_P\to 1/({4\pi^2})$, so that \eqref{eq:CLS_abstract} gives $2C_{\rm LS}\to 2C_{{\rm LS},1}\lambda/(\lambda-1)+{2}/({4\pi^2})$, which is strictly larger than the one-site constant $2C_{{\rm LS},1}$ obtained by independent tensorization: the factor $\lambda/(\lambda-1)$ comes from the Donsker-Varadhan formula, and the term $2C_P$ from Rothaus' lemma \eqref{eq:Rothaus}.

We also observe that  both mechanisms are limited by the same threshold. For a Riesz-type kernel $|x|^{-s}$ on $\T^d$, $\nabla W\in L^{2d/(d+2)}(\T^d,dx)$ if and only if $s<d/2$, and $W\in L^2(\T^d,dx)$ if and only if $s<d/2$. This is the Hilbert-Schmidt regime of \cite{Rosenzweig26}. This is not a coincidence. Indeed, by Parseval's identity, the $H^{-1}(\T^d,dx)$-norm of $\nabla W$ is $$\big(\sum_{k\neq0}|2\pi k\,\hat W(k)|^2/(4\pi^2|k|^2)\big)^{1/2}=\|W\|_{L^2(\T^d,dx)},$$
With respect to the flat measure, the two quantities thus coincide. For the Coulomb kernel $s=d-2$, both conditions hold if and only if $d\le3$, which is why the case $d=3$ is critical for our methods. Beyond this regime, the partition function is not bounded uniformly in $N$ \cite[Proposition~2.2(ii)]{DJ26}, and new ideas are needed.


\section{Appendix}\subsection{Numerical values of the constants}\label{sec.N}

The constants appearing in our results were computed numerically with the help of AI tools (see Remark~\ref{rem:AI}), and we have checked these computations. On $\T^3$, we evaluate $W$ and $\nabla W$ from the Ewald representation \cite{Ewald21} (see also \cite[Section~3]{dLPS80}), and we obtain
$$
W_*=W(\tfrac12,\tfrac12,\tfrac12)\approx-0.802,\qquad\|W\|_{L^1(\T^3,dx)}\approx0.598,\qquad\sup_{|e|=1}\|e\cdot\nabla W\|_{L^{6/5}(\T^3,dx)}\approx4.49 .
$$
On $\T^2$, we evaluate $W$ from its expression in terms of the Jacobi theta function $\vartheta_1$ \cite[Chapter~21]{WhittakerWatson}, and we obtain
$$
W_*=W(\tfrac12,\tfrac12)\approx-0.347,\qquad\|W\|_{L^1(\T^2,dx)}\approx0.263,\qquad\sup_{|e|=1}\|e\cdot\nabla W\|_{L^{4/3}(\T^2,dx)}\approx2.10 .
$$
In both cases, the $L^2$ norm follows from Parseval's identity: 
$$\|W\|_{L^2(\T^3,dx)}=\frac1\pi\big(\sum_{k\neq0}|k|^{-4}\big)^{1/2}\approx1.29 \text{ and } \|W\|_{L^2(\T^2,dx)}=\frac1{2\pi}\big(\sum_{k\neq0}|k|^{-4}\big)^{1/2}\approx0.391.$$

\subsection{On the Sobolev constant $C_{{\rm So}}$} 
\label{sec.Sobo}

\blem\label{le.Sobolev}
For every $w\in H^1(\T^3,dx)$,
$$
\|w\|_{6}\le \|w\|_{2}+\frac2{\sqrt3}\,\|\nabla w\|_{2}
\quad\text{and thus}\quad
\|w\|_{6}^2\le \frac{7}{3}\Big(\|\nabla w\|_{2}^2+\|w\|_{2}^2\Big).
$$
In particular, \eqref{SobolevI} holds with $C_{{\rm So}}\le \frac{7}{3}$.
\nlem

Taking $w\equiv1$ in \eqref{SobolevI} gives $C_{{\rm So}}\geqslant1$, and Lemma~\ref{le.Sobolev} gives $C_{{\rm So}}\le\frac73$.

\bprf 
The proof follows the classical argument of Gagliardo and Nirenberg \cite{Gagliardo58,Nirenberg59}, see e.g.\ \cite[Theorem~9.9]{Brezis11}, adapted to the torus. By density, it suffices to consider $w\in C^\infty(\T^3)$. Identify $\T^3$ with $[0,1)^3$, and for $x\in\T^3$ and $i\in\{1,2,3\}$ write $\hat x_i\in\T^2$ for $x$ with the $i$-th coordinate removed.

\textbf{Step 1.} Let $f\in C^1(\T^3)$, $x\in\T^3$ and $s\in[0,1)$. The points $x$ and $x':=(x_1,\dots,s,\dots,x_3)$, where $s$ is in the $i$-th position, lie on the closed curve $\{x+te_i:\ t\in[0,1)\}$, which is a circle of length $1$ in $\T^3$, and split it into two arcs $\gamma_1$ and $\gamma_2$. Since $f(x)-f(x')=\int_{\gamma_1}\partial_if=-\int_{\gamma_2}\partial_if$, we get
$$
|f(x)-f(x')|\le \min\Big(\int_{\gamma_1}|\partial_if|,\int_{\gamma_2}|\partial_if|\Big)\le \frac12\int_0^1|\partial_if|\,dx_i .
$$
Thus $|f(x)|\le |f(x')|+\frac12\int_0^1|\partial_if|\,dx_i$. Because the left-hand side and the last term do not depend on $s$, integrating this inequality over $s\in[0,1)$ gives
$$
|f(x)|\le  g_i(\hat x_i):=\int_0^1\Big(|f|+\frac12|\partial_if|\Big)(x_1,\dots,t^i,\dots,x_3)\,dt,\qquad i=1,2,3,
$$
where $t^i$ indicates that $t$ is in the $i$-th position. Note that $g_i$ does not depend on $x_i$. As a consequence, $|f(x)|^{3/2}\le \prod_{i=1}^3g_i(\hat x_i)^{1/2}$, and the Cauchy-Schwarz inequality applied successively in $x_1$, $x_2$, $x_3$ (as in the Loomis-Whitney inequality \cite{LoomisWhitney49}) gives
$$
\int_{\T^3}|f|^{3/2}dx\le \prod_{i=1}^3\Big(\int_{\T^2}g_i\,d\hat x_i\Big)^{1/2}
=\prod_{i=1}^3\Big(\|f\|_1+\frac12\|\partial_if\|_1\Big)^{1/2}.
$$
By the arithmetic-geometric mean inequality, it follows that
$$
\|f\|_{3/2}\le \prod_{i=1}^3\Big(\|f\|_1+\frac12\|\partial_if\|_1\Big)^{1/3}\le \|f\|_1+\frac16\sum_{i=1}^3\|\partial_if\|_1 .
$$

\textbf{Step 2.} Apply Step~1 to $f=w^4$. Then, $\|f\|_{3/2}=\|w\|_6^4$. By H\"older's inequality, $\|f\|_1\le \|w\|_6^3\|w\|_2$ and $\|\partial_if\|_1=4\int|w|^3|\partial_iw|\,dx\le 4\|w\|_6^3\|\partial_iw\|_2$. Since $\sum_{i=1}^3\|\partial_iw\|_2\le \sqrt3\,\|\nabla w\|_2$, we get that:
$\|w\|_6^4\le \|w\|_6^3(\|w\|_2+ 2 \|\nabla w\|_2/{\sqrt3} )$,
which is the first inequality. The second one follows from $(a+\beta b)^2\le (1+\beta^2)(a^2+b^2)$ with $\beta= 2/{\sqrt3}$.
\nprf
%


\subsection{Proof of Lemma \ref{le.HJ}}
\label{sec.HJ2}

\bprf(Lemma \ref{le.HJ})
Let $\overline f$ be the $\zz^3$-periodic lift of $f$ to $\rr^3$, which is bounded and Lipschitz with the same constant. Since $d(x,y)=\min_{k\in\zz^3}|x-y+k|$, the lift of $Q_tf$ is
$$
\overline{Q_tf}(x)=\min_{y\in\rr^3}\Big\{t\,L\Big(\frac{x-y}{t}\Big)+\overline f(y)\Big\},\qquad L(q):=|q|^2,\quad x\in\rr^3,
$$
which is the Hopf-Lax formula \cite[Section~3.3.2, Theorem~4]{Evans10} with the Lipschitz initial datum $\overline f$ and the Hamiltonian
$H(p)=L^*(p)=\sup_{q\in\rr^3}\{p\cdot q-L(q)\}=\frac14|p|^2$, see \cite[Section~3.3.2, Theorem~3]{Evans10}. Because $H$ is smooth, convex and superlinear, the results of \cite[Sections~3.3.2 and~10.3.4]{Evans10} apply.

\emph{(a)} The function $t\mapsto f(y)+\frac1td^2(x,y)$ is nonincreasing for every $y$, and so is its infimum over $y$. Taking $y=x$ gives $Q_tf\le f$, and $Q_tf\ge\min f$ is clear.

\emph{(b)} By \cite[Section~3.3.2, Lemma~2]{Evans10}, $(t,x)\mapsto Q_tf(x)$ is Lipschitz on $[0,\infty)\times\T^3$ and $Q_0f=f$. Moreover $\|Q_tf\|_{\rm Lip}\le\|f\|_{\rm Lip}$ by \cite[Section~3.3.2, (28) in the proof of Lemma~2]{Evans10}.

\emph{(c)} By \cite[Section~3.3.2, Theorem~5]{Evans10}, $Q_tf$ solves \eqref{lem51e} at every point of differentiability of $(t,x)\mapsto Q_tf(x)$. By \emph{(b)} and Rademacher's theorem, these points have full measure (see also \cite[Section~3.3.2, Theorem~6]{Evans10}), and Fubini's theorem gives the last assertion.

\emph{(d)} By \cite[Theorem~1.3]{Villani21}, $W_2^2(\nu,\mu)$ is the supremum of $\int \psi\,d\nu-\int f\,d\mu$ over pairs of continuous functions such that $\psi(y)-f(x)\le d^2(x,y)$ for all $x,y$. For a given $f$, the largest admissible $\psi$ is $Q_1f$, and the supremum may be taken over $f\in C(\T^3)$ with $\psi=Q_1f$. Finally, $|Q_1f-Q_1g|\le\|f-g\|_\infty$, and by uniform approximation the supremum may be restricted to $f\in C^\infty(\T^3)$.

\emph{(e)} Since $\overline f$ is bounded and Lipschitz, the proof of \cite[Section~10.3.4, Theorem~3]{Evans10}, which only uses the convexity of $H$, shows that $\overline{Q_tf}$ is a bounded and Lipschitz viscosity solution of $u_t+ |\nabla u|^2/4=0$ on $\rr^3\times(0,\infty)$ with $u(\cdot,0)=\overline f$. Uniqueness among bounded and uniformly continuous viscosity solutions follows from \cite[Theorem~4.1]{CEL84}, which only requires $H$ to be continuous (whereas the uniqueness theorem \cite[Section~10.2, Theorem~1]{Evans10} assumes $H$ to be globally Lipschitz). Finally, since $\overline{Q_tf}$ is $\zz^3$-periodic and the viscosity property is local, the same holds on $\T^3$, by lifting solutions to periodic ones on $\rr^3$.
\nprf

\subsection{A uniform $W_1H$ inequality in every dimension}\label{sec.W1H}

Let $d\geqslant2$ and $D_d:={\rm diam}(\T^d)=\frac{\sqrt d}{2}$. We equip $(\T^d)^N$ with the metric $d_{\ell^1}(x,y):=\sum_{i=1}^Nd(x^i,y^i)$.

\bprop\label{P_W1H}
Let $W$ be an even function on $\T^d$ with $\int_{\T^d}W\,dx=0$, $W\geqslant W_*>-\infty$ and $\nabla W\in L^1(\T^d,dx)$, and let $\mu_N$ be given by \eqref{Coulomb-measure} on $(\T^d)^N$. Assume that
$$
r_{W_1}:=\frac{D_d}{2}\,\Gamma\,e^{-\Gamma W_*}\,\|\nabla W\|_{L^1(\T^d,dx)}<1 .
$$
Then, for every $N\geqslant2$ and every probability measure $\nu$ on $(\T^d)^N$,
\bbeq\label{eq:W1H}
W_{1,d_{\ell^1}}(\nu,\mu_N)^2\le\frac{D_d^2\,N}{2(1-r_{W_1})^2}\,H(\nu|\mu_N).
\nneq
As a consequence, for every $1$-Lipschitz function $f$ on $\T^d$ and every $t>0$,
$$
\mu_N\Big(\Big|\frac1N\sum_{i=1}^Nf(x^i)-\ee_{\mu_N}\Big[\frac1N\sum_{i=1}^Nf(x^i)\Big]\Big|\geqslant t\Big)\le2\exp\Big(-\frac{2(1-r_{W_1})^2Nt^2}{D_d^2}\Big).
$$
\nprop

The assumptions hold for the Coulomb kernel on $\T^d$ for every $d\geqslant2$, since $|\nabla W(x)|\le C|x|^{1-d}$ near $0$, and for the Riesz kernels $|x|^{-s}$ with $s<d-1$.
Chafa\"i, Hardy and Ma\"ida \cite{CHM18} obtained concentration inequalities for the empirical measure of Coulomb gases in every dimension $d\geqslant2$, for an interaction $\frac\beta2\sum_{i\neq j}g(x^i-x^j)$ with $\beta>0$ fixed, where the energy dominates the entropy. In the mean-field regime $\beta=\frac{\Gamma}{N-1}$ considered here, their bound becomes trivial, and Proposition~\ref{P_W1H} provides a Gaussian concentration at small charge.

\bprf
Let $i\neq j\in \{1,\ldots,N\}$, and let $x^{-i},y^{-i}$ differ only at the coordinate $j$. Let $(\mu_t)_{t\in[0,1]}$ be defined as in Step~0 of the proof of Proposition~\ref{P-coef_W2}, with $\T^3$ replaced by $\T^d$. Step~1 of that proof only uses $\rho_\tau\le e^{-\Gamma W_*}$ and $\nabla W\in L^1(\T^d,dx)$, and it holds on $\T^d$. If $f$ is $1$-Lipschitz, then $\osc(f)\le D_d$, and Step~1{(i)} of that proof applied to $f/D_d$ gives
$$
\mu_i(f|y^{-i})-\mu_i(f|x^{-i})\le\frac{\Gamma D_d}{2(N-1)}\,e^{-\Gamma W_*}\,\|\nabla W\|_{L^1(\T^d,dx)}\,d(x^j,y^j).
$$
Recall that (see Section~\ref{sec:abstractW})
$$
c^{d,1}_{ij}=\sup_{x^{-i}=y^{-i}\text{ off }j}\frac{W_{1,d}\big(\mu_i(\cdot|x^{-i}),\mu_i(\cdot|y^{-i})\big)}{d(x^j,y^j)},
\qquad
W_{1,d}(\mu,\nu)=\sup_{\|f\|_{\rm Lip}\le1}\big(\mu(f)-\nu(f)\big),
$$
where the second identity is the Kantorovich-Rubinstein duality. Taking the supremum over all $1$-Lipschitz functions $f$ in the previous inequality, we get that $c^{d,1}_{ij}\le\frac{r_{W_1}}{N-1}$, and thus $\sup_j\sum_{i\neq j}c^{d,1}_{ij}\le r_{W_1}<1$. 
Moreover, by Pinsker's inequality, every probability measure $\nu$ on $\T^d$ satisfies $W_{1,d}(\nu,\mu)\le D_d\|\nu-\mu\|_{\rm tv}\le D_d\sqrt{H(\nu|\mu)/2}$, for every probability measure $\mu$ on $\T^d$. Thus, for every probability measure $\mu$ on $\T^d$ (in particular for every conditional law of $\mu_N$ given some of the coordinates), we have
$$
W_{1,d}(\nu,\mu)^2\le 2K\,H(\nu|\mu),\qquad K:=\frac{D_d^2}{4},
$$
for every probability measure $\nu$ on $\T^d$. This is the one-site $W_1H$ inequality required in \cite[Theorem~4.3]{Wu_AOP06}, which holds here for free since $\T^d$ is bounded, uniformly in $N$ and in the positions of the poles. 
Because the matrix $C=(c^{d,1}_{ij})$ satisfies $\|C\|_1:=\sup_j\sum_{i\neq j}c^{d,1}_{ij}\le r_{W_1}<1$, we can use  \cite[Theorem~4.3]{Wu_AOP06} which  implies \eqref{eq:W1H}, together with 
$$\ee_{\mu_N}[e^{\lambda(F-\ee_{\mu_N}[F])}]\le\exp\Big(\frac{\lambda^2D_d^2N\alpha^2}{8(1-r_{W_1})^2}\Big), \quad \lambda\in\rr,$$ for every $F$ which is $\alpha$-Lipschitz in each coordinate. Taking $F=\frac1N\sum_if(x^i)$, i.e.\ $\alpha=\frac1N$, and optimizing Chebyshev's exponential inequality in $\lambda$ gives the concentration bound.
\nprf

\begin{remark}
Inequality \eqref{eq:W1H} is a transport-entropy inequality of Talagrand type with cost $W_1$ (see \cite{DGW04,Wu_AOP06}). Since every probability measure on a bounded space satisfies such an inequality by Pinsker's inequality, the content of Proposition~\ref{P_W1H} is the order $N$ of the constant, which is the same as for product measures. Inequality \eqref{eq:W1H} does not imply the Poincar\'e inequality, but by \cite{BobkovGotze99} it is equivalent to a Gaussian concentration of Lipschitz functions of the configuration, whereas the Poincar\'e inequality only gives an exponential concentration.
\end{remark}

\vskip 12pt
\noindent
{\small
\textbf{Acknowledgement}. A.G. is supported by the ANR-23-CE-40003, Conviviality, and has benefited from a government grant managed by the Agence Nationale de la Recherche under the France 2030 investment plan ANR-23-EXMA-0001. W.L. is supported by NSFC 12471143. B.N. is supported by the ANR-25-CE40-6875.}

\bibliographystyle{plain}
\bibliography{Coulomb_propagation}

\end{document}